\documentclass[12pt]{article}

\usepackage{latexsym,amssymb,amsgen,amsmath,amsxtra, amsfonts, amsthm,amscd}
\usepackage{pstricks,pst-tree,pstricks-add}

\usepackage{graphics}

\usepackage{graphicx}

\usepackage{subfigure}

\usepackage{authblk}

\usepackage{color}

\newtheorem{theorem}{Theorem}

\newtheorem{corollary}[theorem]{Corollary}

\newtheorem{lemma}[theorem]{Lemma}

\newtheorem{remark}[theorem]{Remark}

\newtheorem{defin}[theorem]{Definition}

\newtheorem*{thma}{Theorem (Hutchinson)}
\newtheorem*{thmb}{Theorem (Elton)}

\newtheorem*{ass1}{Assumption 1}
\newtheorem*{ass2}{Assumption 2}

\def\be{\begin{equation}}
\def\ee{\end{equation}}
\def\ba{\begin{array}}
\def\ea{\end{array}}

\def\d{\hbox{\sl dist}}

\def\Proof{\noindent{\bf Proof.\ \ }}

\begin{document}

\title{Dynamics with Choice}
\author{
Lev Kapitanski \thanks{
The work was supported in part by the
NIH  Director's Initiative,
1 P20 RR020770-01.}
}
\author{Sanja \v Zivanovi\'c}
\affil{Department of Mathematics\\
University of Miami\\
Coral Gables, FL 33124, USA}

\date{\today}


%

%

%
\maketitle


\begin{abstract}
Dynamics with choice is a generalization of discrete-time dynamics where instead of the same evolution operator at every time step there is a choice of operators to transform the current state of the system. Many real life processes studied in chemical physics, engineering, biology and medicine, from autocatalytic reaction systems to switched systems to cellular biochemical processes to malaria transmission in urban environments, exhibit the properties described by dynamics with choice. We study the long-term behavior in dynamics with choice. We prove very general results on the existence and properties of
global compact attractors in dynamics with choice.
In addition, we study the dynamics with restricted
choice when the
allowed sequences of operators correspond to
subshifts of the full shift.
One of practical consequences of our results
is that when the parameters of a discrete-time  system are not known exactly and/or are subject to change due to internal instability, or a strategy,
or Nature's intervention, the long term behavior
of the system may not be correctly described by a system with ``averaged" values for
the parameters. There may be a Gestalt effect.
\end{abstract}


\section{Introduction}\label{intro}

\begin{figure}[hb!]
\centering
\includegraphics[width=.7\textwidth]{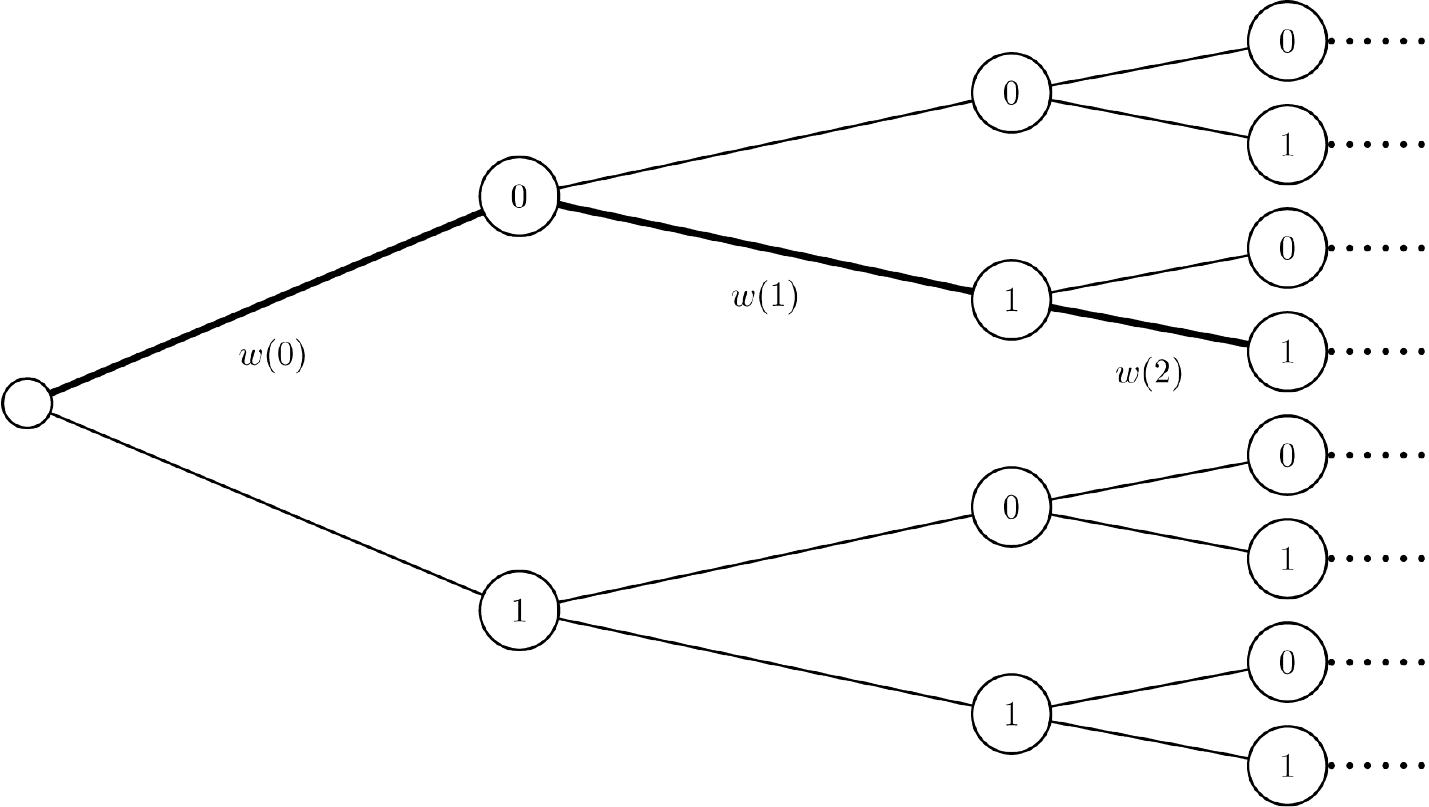}
\caption{The tree of choices in the case of two operators}
\label{f1}
\end{figure}

Mathematical setting for discrete dynamics is a
space $X$ and a map $S:\,X\to X$. The space $X$ is
the state space, the space of all possible states of the system. The map $S$, the evolution operator,
defines the change of states over one time step:
$x\in X$ at time $t = 0$ evolves into $S(x)$ at $t = 1$,
$S(S(x))$ at $t = 2$, \dots, $S^{n}(x)$ at $t = n$, etc.
If instead of one operator, $S$, we have a choice of evolution operators, $S_0$, $S_1$, \dots,
$S_{N-1}$, and at every time step choose one of them,
then we have a {\it dynamics with choice}.
One way to visualize the multitude of choice
through time is to generate the infinite tree of choices.
This is an infinite rooted tree in which the root has
$N$ children, every child has $N$ children, and so on.
The root corresponds to $t=0$, its children correspond to
$t=1$, the children of the children correspond to $t=2$, etc.
At every step, the children of each node
are labeled $0$ through $N-1$.
Beginning at the root infinite branches (paths, strategies) represent
the possible choices: for example, in Figure \ref{f1}
we choose the path $w$ that starts with 011... (bold edges).
For this choice,
the first few points in the trajectory of a point $x_0\in X$ are
$x_1 = S_0(x_0)$, $x_2 = S_1(x_1) = S_1(S_0(x_0))$,
$x_3 = S_1(x_2) = S_1(S_1(S_0(x_0)))$, etc.
It is natural to encode the infinite paths
(beginning at the root) by one-sided infinite words
(strings, sequences)
on $N$ symbols. If $w$ is such sequence, it is
convenient to align it with the set of non-negative integers
${\mathbb Z}_{\ge 0}$
and denote by $w(k)$ the $(k+1)$-st letter of $w$, i.e.,
$w = w(0)w(1)w(2)\ldots$. Thus, $w(0) = 0$, $w(1) = 1$, $w(2) = 1$,
are the first three symbols of the path
$w = 011\ldots$.

\noindent$\bullet$\ \ In this paper we study dynamics with choice, i.e., the dynamics of
points and subsets of $X$ along all possible paths simultaneously. We will explain what this means momentarily.
Here we would like to emphasize that, from the point of view of long-term behavior, dynamics with choice,
in general, is not the same as the union of trajectories along
different infinite paths. We will return to this point later
when we talk about the Gestalt effect.
\bigskip

Let $\Sigma$ be the one-sided shift on $N$ symbols, \cite{Kitchens}. This means that, first,
$\Sigma$ is the set of all one-sided infinite
strings $w = w(0)w(1)w(2)\dots $, where each $w(j)$
is a symbol from the list of $N$ symbols $[0, 1, \dots, N-1]$,
and second,
there is the shift operator $\sigma:\,\Sigma\to \Sigma$
acting by erasing the first symbol,
$\sigma(w) = w(1)w(2)w(3)\dots$. Given the state space $X$
and operators $S_0, \dots, S_{N-1}$,
we define the corresponding {\it dynamics with choice}
as the discrete dynamics
on the state space the product ${\frak X} = X\times \Sigma$
with the following evolution operator:
\be\label{rule}
{\frak S}:\;(x, w)\mapsto (S_{w(0)}(x),\,\sigma(w))\,.
\ee
One can think of a $w\in\Sigma$ as a plan, a strategy, or as Nature's intervention. Dynamics with choice is a language  to  describe  processes where different strategies could be applied or happen.
Most of mathematical models in natural sciences
and engineering
are expressed in terms of differential equations.
Those equations are often continuous limits of
discrete equations. Continuous case is easier for
qualitative analysis. However, there are situations
where discrete equations describe the processes better. Every realistic model comes with
parameters. We are interested in situations where
parameters may change
due to, e.g.,
internal instability or outside intervention.
In an illustrative example in section
\ref{ex}, the coefficients $a$ and $b$
are proportional to the biting rate of mosquitoes
which depends, for example, on temperature
and humidity which may change from day to day
and during the day.

\medskip

In this paper we study long-term regimes in dynamics with choice. More specifically, we define and
study global
compact attractors in dynamics with choice.
By a global compact attractor we mean {\bf the minimal compact set that attracts all bounded sets},
see section \ref{attr} for definitions and references.
Thinking in terms of a model with parameters,
assume we know that for each admissible fixed
(in time) set
of parameters the system possesses a global compact
attractor. What happens when the parameters
switch between admissible values? Is there an attractor? How is it related to attractors corresponding to fixed parameters?
Is there a Gestalt effect?
These are the questions we address in this paper.

\medskip

There are many real life and engineered
systems that switch between different modes
of operation (the so-called hybrid systems).
When the behavior in each mode is modeled
using continuous dynamics and the transitions
are viewed as discrete-time events, such systems
are called switching or switched. Analysis and
especially control of switching systems is
an area of intensive research, see, e.g.,
Liberzon's book \cite{Lib} and the survey
by Margaliot \cite{Marg}. There is a natural
affinity between switching systems and dynamics with choice (see, e.g., \cite{Kloed}), but we will not
explore it at this time.

\bigskip

Readers familiar with iterated function systems,
\cite{Hutchinson,Barnsley-1},
may wonder if there is a connection between
iterated function systems and  dynamics with choice.
Indeed there is, but we have to establish it
(in section \ref{ifs}).

A (general) Iterated Function System
(IFS)\footnote{We use the abbreviation IFS
for single and IFSs for plural forms.}
can be viewed as a discrete dynamics on
the space $2^X$ (of subsets of $X$). The operators  $S_0$, $S_1$, \dots,
$S_{N-1}$ define evolution on $2^X$ by means of the Hutchinson-Barnsley operator:
\be\label{HB}
F\,:\;A\mapsto F(A) :=
S_0(A)\cup S_1(A)\cup\cdots\cup S_{N-1}(A)\,.
\ee
Following a long-standing tradition,
people studying dynamics are first of all interested
in fixed points. In the case of an IFS, those are
the fixed points of the Hutchinson-Barnsley operator.
As has been well illustrated by Barnsley,
for many simple IFSs on the plane one can (use a computer to)
plot their compact fixed points (sets) and
obtain fascinating fractals, see \cite{Barnsley-1,Barnsley-2}. Generating  fractals
is one of the main motivations in the study of
IFSs. In some papers a fractal is defined as the
compact invariant set of an IFS, see \cite{Barnsley-2} and references therein.

To prove that an IFS does have a fixed point,
the general definition should be made more specific.
One needs to specify the properties of
the space $X$;
the space $2^X$ should be narrowed to an appropriate
class of subsets; assumptions should be made on
the operators $S_0$, $S_1$, \dots,
$S_{N-1}$. As an example we state the
original result of Hutchinson, \cite[Section 3]{Hutchinson}.
\begin{thma}
Let $X$ be a complete metric space
(with metric $d$). Denote by
$\bar{\cal B}(X)$ the space of all non-empty closed bounded
subsets of $X$. Assume that each operator
$S_0$, $S_1$, \dots,
$S_{N-1}$ is a strict contraction (i.e.,
there is a number  $\gamma\in (0,1)$ such that
$d(S_j(x),\,S_j(y))\le \gamma\,d(x,\,y)$ for every pair $x, y\in X$ and for all $j$). Define the
evolution operator
$\bar F:\,\bar{\cal B}(X)\to \bar{\cal B}(X)$ by the formula
\be\label{barF}
\bar F : \,A\mapsto \bar F(A) =
\overline{S_0(A)\cup S_1(A)\cup\cdots\cup S_{N-1}(A)}\,.
\ee
Then there exists a unique fixed point
$K\in \bar{\cal B}(X)$ of $\bar F$.
Viewed as a subset of $X$, the set $K$ is compact.
Also, $K$ attracts every closed bounded subset of $X$ in the sense that, for any  $C\in \bar{\cal B}(X)$,
\[
d_H({\bar F}^n(C),\,K) \to 0\qquad\hbox{as}\quad
n\to\infty\,,
\]
where $d_H$ is the Hausdorff distance.
\end{thma}
The IFS with contractive operators $S_j$
are called hyperbolic.
Over the years this result has been generalized
in many different
directions (different assumptions on $X$ and/or $S_j$), see \cite{AFGL} for references.

\bigskip

The iterated function systems with probabilities
and the `chaos games' in particular
show an apparent link to dynamics with choice.
Recall that {\it an iterated function system with probabilities} is an IFS $(X;\,S_0,\dots, S_{N-1})$
together with probabilities $p_0,\,p_1, \dots,\,p_{N-1}$
assigned to the operators $S_0$, $S_1$, \dots,
$S_{N-1}$, where each $p_j > 0$ and $p_0 + p_1 + \dots + p_{N-1} = 1$, \cite{Barnsley-1}.
The random iteration algorithm (aka the chaos game,
\cite{Barnsley-2})
starts with the choice of initial
state $x_0\in X$. Next, define recursively
$x_{n+1}$ by choosing its value from the set
$\{S_0(x_n),\,S_1(x_n),\dots, S_{N-1}(x_n)\}$
with respective probabilities
$p_0,\,p_1, \dots,\,p_{N-1}$.
The choice of operators thus will be encoded in
some strategy
$w = w(0)w(1)w(2)\dots\in \Sigma$, i.e.,
$x_{n+1} = S_{w(n)}(x_n)$.  To show that
the sequence $(x_n)$ is determined by $x_0$ and $w\in \Sigma$, we  write
$(x_n(x_0,\,w))$.  Consider the averages
of the delta-measures  concentrated at the points
$x_n$.
It turns out that (under certain conditions) the averages $n^{-1}\left(\delta_{x_0} + \delta_{x_1} + \dots + \delta_{x_{n-1}}\right)$ converge weakly
to the invariant measure of the IFS.
More precisely,
consider the following Markov operator
on the space of probability measures on $X$:
\[
M :\,\nu\mapsto M(\nu) = \sum_{j=0}^{N-1} p_j\,S_j(\nu)\,,
\]
where $S_j(\nu)(A) = \nu(S_j^{-1}(A))$ for measurable sets $A\subset X$. Hutchinson showed in
\cite{Hutchinson} that under the assumptions of the Hutchinson Theorem  there exists a unique fixed point  $\mu$ of the Markov operator and the support of $\mu$ is the fixed point $K$ of the IFS
$(X;\,S_0,\dots, S_{N-1})$. The following theorem
was proved
by Elton, \cite{Elton}, with later simplifications by
Forte and Mendivil, \cite{FM}.
\begin{thmb}
Assume $X$ is a compact metric space, the operators
$S_0,\dots, S_{N-1}$ are strict contractions on $X$,
and $p_0,\,p_1, \dots,\,p_{N-1}$ are the probabilities. Let $\mu$ be the corresponding invariant probability measure.
Then, for any continuous function
$f:\,X\to \mathbb R$ and any $x_0\in X$,
\[
\lim\limits_{n\to\infty}\;\;
\frac1{n}\;\sum\limits_{i = 0}^{n-1}
f(x_i(x_0,\,w))\,=\,\int_X f(x)\,d\mu(x)\,,
\]
for almost all strategies $w\in\Sigma$ with respect to the product probability measure on
$\Sigma = [0, 1,\dots, N-1]^{\mathbb N} $ induced
by the distribution $[p_0, p_1,\dots, p_{N-1}]$
on each factor.
\end{thmb}

Returning to dynamics with choice, we repeat that
our interest has not been motivated by fractals.
We would like to understand the long-term behavior
in dynamics with choice.
We assume that $X$ is a complete
metric space (with metric $d$), the operators
$S_0, \dots, S_{N-1}$ are continuous,
and each of the
(semi)dynamical systems $(X, d,\,S_j)$ possesses a
global compact attractor.
Consider the corresponding dynamics with choice
as the dynamics on the product metric space
\footnote{$\Sigma$ can be equipped with a metric
making it a compact metric space, see Section
\ref{sigma} for a specific choice. We denote here by
$\hbox{\sl dist}$ the corresponding product-metric on $X\times \Sigma$. }
${\frak X} = X\times \Sigma$ generated by the operator ${\frak S}$ acting according to the rule  (\ref{rule}). From general theory (see section \ref{attr})
we know that a system ought to enjoy certain
compactness and dissipativity properties
in order for it to possess the global compact attractor.


In general, even when the individual systems $(X, d,\,S_j)$ do have attractors, the  system
$({\frak X}, \hbox{\sl dist}, {\frak S})$
will not have a global compact attractor. There are
several reasons why.
One counter-example we borrow from  \cite{AF-2} (where it is used in
the context of IFS). Take $X = \mathbb R$ with standard metric $d$
and define two maps, $S_0$ and $S_1$, as follows:
\[
S_0(x) = \left\{
\begin{array}{ll}
0, & \textrm{if $x \le 0$},\\
- 2 x, & \textrm{if $x > 0$}
\end{array}
\right.
\qquad
S_1(x) = \left\{
\begin{array}{ll}
- 2x, & \textrm{if $x \le 0$},\\
0, & \textrm{if $x > 0$}
\end{array}
\right.
\]
Each of the systems $(X, d, S_j)$ has the global compact attractor,
a singleton $\{0\}$. At the same time, the trajectory $x_n = S_{w(n-1)}\circ S_{w(n-2)}\circ\dots \circ S_{w(0)}(x_0)$ corresponding to the periodic string $w = 010101\dots$ is unbounded
for any initial point $x_0\neq 0$. Hence, there is no compact attractor attracting $(x_0, w)$.

The second example is infinite-dimensional.
Let $B_0=B_0(p_0)$ and $B_1=B_1(p_1)$ be two  disjoint closed unit balls
centered at
$p_0$ and $p_1$
in an
infinite-dimensional Banach space. Let $X = B_0\cup B_1$. Define the maps $S_0$ and $S_1$ as follows: on $B_0$ the map $S_0$ is a contraction and it maps $B_1$ to $B_0$;
the map $S_1$ is a contraction on $B_1$ and  maps $B_0$ to $B_1$:
\[
S_0(x) = \left\{
\begin{array}{ll}
p_0 + \frac12\,(x - p_0), & \textrm{if $x \in B_0$},\\
p_0 + (x - p_1), & \textrm{if $x \in B_1$}
\end{array}
\right.
\qquad
S_1(x) = \left\{
\begin{array}{ll}
p_1 + \frac12\,(x - p_1), & \textrm{if $x \in B_1$},\\
p_1 + (x - p_0), & \textrm{if $x \in B_0$}
\end{array}
\right.
\]
The system $(X, d, S_0)$ does have the global compact attractor, $\{p_0\}$, and
$(X, d, S_1)$ does have the global compact attractor, $\{p_1\}$. The corresponding dynamics with choice, $({\frak X}, \hbox{\sl dist}, {\frak S})$, does have the global {\it closed} attractor,
namely, ${\frak X}$, but does not have
the global {\it compact} attractor.

In the first example, the maps are compact (which is good), but they do not have a joint bounded absorbing set
(lack of dissipativity in $({\frak X}, \hbox{\sl dist}, {\frak S})$). In the second example,
there is a joint bounded absorbing set,
$B_0\cup B_1$, but there is not enough
compactness (the maps $S_j$ are not compact, not
contracting,
and, more generally, not condensing).

These examples show what kind of situations do not allow
global compact attractors in the dynamics with choice.
Thus, we make additional assumptions. First, we assume that
there exists a bounded absorbing set that absorbs every bounded set regardless of the strategy.
In applications, absorbing set is usually a ball of the radius that depends on the parameters of the model. Our ``dissipativity" assumption means that there is a common  estimate on the radius for different values of
the parameters.


\begin{ass1}
There is a closed, bounded set ${\bf B}\subset X$ such that for every bounded $A\subset X$ there exists $m(A)>0$ such that
$S_{w(n-1)}\circ S_{w(n-2)}\circ \cdots \circ
S_{w(0)}\,(A)  \,\subset {\bf B}$
for every word $w = w(0)w(1)\dots w(n-1)$ of length $n\ge m(A)$.
\end{ass1}

Our second, ``compactness" assumption is that each of the operators $S_j$ is condensing with respect to a common measure of noncompactness. This assumption covers practically all  situations
encountered in applications: contractions, compact operators, and compact plus contractions.
As their name suggests,
measures of noncompactness measure how far a set is from being compact. There are several
different measures of noncompactness in use,
\cite{AKPRS}. For example,
the Hausdorff measure of noncompactness of
a set $A$ is the infimum
of $\epsilon > 0$ such that $A$ has a finite
$\epsilon$-net.  In this paper we use only very general properties
shared by
all popular measures of noncompactness, see
Definition \ref{nc} in section \ref{mnc} below.

Let $\psi$ be a measure of noncompactness (as in
Definition \ref{nc}). An operator $S:\,X\to X$
is condensing with respect to $\psi$
iff $\psi(S(A)) < \psi(A)$ for any non-compact set
$A$, and $\psi(S(A)) = \psi(A) = 0$ if $A$ is compact. Our second general assumption is this.
\begin{ass2}
Each operator
$S_j$ is $\psi$-condensing.
\end{ass2}

\medskip

In section \ref{dwc} we prove the following result.


\begin{theorem}\label{one}
Let $X$ be a complete metric space
and let
$S_0, S_1,\dots, S_{N-1}$ be continuous, bounded
(i.e., take bounded sets to bounded sets) maps
$X\to X$. In addition, let assumptions 1 and 2 be
satisfied.
Then the system $({\frak X}, \hbox{\sl dist},\,{\frak S})$ has a global compact attractor, ${\frak M}$. (That $\frak M$ is the global compact attractor
means that $\frak M$ is the smallest compact in
$\frak X$
attracting every bounded set in $\frak X$.)
\medskip

The attractor $\frak M$ has the following properties.
\begin{enumerate}
\item[(1)] $\frak M$ is (strictly) invariant:
${\frak S} ({\frak M}) = {\frak M}$.
\item[(2)] $\frak M$ is the union of all closed bounded sets $A\subset{\frak X}$ with the property $A\subset {\frak S}(A)$.
\item[(3)] $\frak M$ is the maximal closed set
with the property $A\subset {\frak S}(A)$; in particular,
$\frak M$ is the maximal (strictly) invariant closed  set.
\item[(4)] Through every point $(x, w)\in {\frak M}$
passes a complete trajectory. This means there exists
a two-sided sequence $\dots, x_{-2}, x_{-1}, x_0, x_1, x_2,\dots $ of points in $X$ and a two-sided
infinite string $\dots s(-2)s(-1)s(0)s(1)s(2)\dots$
such that $x(0) = x$ and
$s(0)s(1)s(2)\dots = w(0)w(1)w(2)\dots$  and such that
$
S_{s(n)}(x_n) = x_{n+1}
$
for every integer $n$.
\item[(5)] $\frak M$ is the union of all complete, bounded trajectories in ${\frak X}$.
\end{enumerate}
\end{theorem}

Given the state space $X$
and operators $S_0,\dots, S_{N-1}$, there are
two ways of describing dynamics generated by
the corresponding IFS.
First, one can follow the trajectories of bounded
subsets of $X$ under the iterations of the
Hutchinson-Barnsley map $\bar F$, see
(\ref{barF}). We denote such system by
$(X, d, {\bar F})$. The notion of the global
compact attractor as the minimal compact set that
attracts all bounded sets, is well-defined for
$(X, d, {\bar F})$. The second possibility is
to choose the space of  closed bounded sets,
$\bar{\cal B}(X)$, as the state space of the system
and study the dynamics of its points
under the iterations of $\bar F$. As a rule, $\bar{\cal B}(X)$ is equipped with
the Hausdorff distance $d_H$.
Thus we obtain the second system,
$(\bar{\cal B}(X), d_H, {\bar F})$. It turns out that
from the point of view of global compact attractors
the dynamical system $(\bar{\cal B}(X), d_H, {\bar F})$
is not very interesting (because convergence in
the Hausdorff metric is too strong).
It possesses an attractor (in the sense
we use here) essentially only if  the maps $S_j$
are contractions, so then the attractor is just
one point in  $\bar{\cal B}(X)$. For more general
$S_j$, it makes more sense to study the fixed points
of $\bar F$.


 In sections \ref{ias} and \ref{inter} we establish the following
connection between the dynamics with choice and
the corresponding IFS.
\begin{theorem}\label{two}
Make the same assumptions on the space $X$
and operators $S_0,\dots, S_{N-1}$ as in Theorem
\ref{one}. Then
\begin{enumerate}
\item[(1)] The IFS\ \  $(X, d, {\bar F})$
does have a global compact attractor, $K$.
\item[(2)] The  set $K $
is the largest compact set in $X$ which is invariant
under the Hutchinson-Barnsley map $\bar F$,
$K = {\bar F} (K)$.
\item[(3)] The attractor $\frak M$ of the dynamics with choice has the following product structure:
\[
\frak M = K\times \Sigma\,.
\]
\end{enumerate}
\end{theorem}

\noindent In the extensive literature on IFSs the main question is the existence of ``the fractal", i.e. the maximal compact
set invariant under the Hutchinson-Barnsley operator
$\bar F$. This corresponds to the second assertion of
our Theorem \ref{two}. We believe that viewing
``the fractal" of an IFS as the attractor of
the dynamical system
$(X, d, {\bar F})$ is beneficial to the theory of
IFSs. This approach, in particular, points to the
``right" assumptions on the space $X$ and the operators $S_j$.

Iterated function systems with compact
(possibly multi-valued) operators have been  considered previously, see, e.g., \cite{AF-2}. The statement of
Theorem 5.8 in \cite{AF-2} which establishes the existence of a compact set invariant under
$\bar F$, needs some additional (dissipativity) assumption such as our Assumption 1, for example.
The IFS with  condensing (and multi-valued, in addition) operators
have been considered by Le\'sniak, \cite{Lesniak}
and Andres et al., \cite{AFGL}.
The assumptions of Theorem 3 in \cite{AFGL} require
that the image ${\bar F}(X)$ of the whole space $X$ be bounded. The word ``minimal" referring to the
``fractal" in \cite[Theorem 3]{AFGL} should probably be replaced by ``maximal", see also
\cite[Theorem 3]{Lesniak}.

\bigskip

 Our assumptions on the state space and the operators guarantee that, for every fixed $j = 0, \dots, N-1$, the discrete dynamics generated on $X$ by $S_j$ does possess the global compact attractor (in $X$). More generally, as we show in
sections \ref{ias} and \ref{inter}, it makes sense
to define {\it individual attractors},
${\cal A}_w$ corresponding to every string
(infinite path in the tree of choices) $w\in \Sigma$. The attractors generated by each $S_j$ correspond
to ``constant" strings, $w = jjj\dots$.
It is not hard to see that such attractors do not exhaust the attractor (fractal) $K$.
There are situations when the union of all
${\cal A}_w$ is $K$
(this happens, in particular,
when $S_j$'s are strict contractions).
However, in general, the union
$\bigcup\limits_{w\in \Sigma}{\cal A}_w$
is {\it strictly smaller} than $K$.
We give an example of this in section \ref{inter}.
In the cases when
$\bigcup\limits_{w\in \Sigma}{\cal A}_w$
is  strictly smaller than $K$ we say that there is
{\bf a Gestalt effect}, i.e., ``the whole is greater than the sum of its parts." This is a new phenomenon. It has not been observed
in the framework of Iterated Function Systems because, as we show
in Lemma \ref{equality}, the Gestalt effect
cannot occur when operators $S_j$ are contractions.


\bigskip

An important generalization of dynamics with choice
is {\it dynamics with restricted choice}.
The name should indicate that not all strategies
(sequences $w = w(0)w(1)\dots\in\Sigma$) are allowed. In particular, we
consider the sets in $\Sigma$ that are closed and
shift invariant, i.e., {\it subshifts}, see \cite{Lind-Marcus,Kitchens}. Given a subshift
$\Lambda\subset\Sigma$, we consider the dynamics on the product-space ${\frak X}_\Lambda = X\times \Lambda$ generated by the map ${\frak S}$ as in
(\ref{rule}).


\begin{theorem}\label{three}
Let the space $X$ and the operators $S_0$,\dots, $S_{N-1}$ satisfy  Assumptions 1 and 2. Let $\Lambda$
be a (one-sided) subshift of $\Sigma$. Consider the
dynamical system $({\frak X}_\Lambda, \hbox{\sl dist}, {\frak S})$.
\begin{enumerate}
\item[(1)] The dynamical system $({\frak X}_\Lambda, \hbox{\sl dist}, {\frak S})$ does possess a global compact attractor, ${\frak M}_\Lambda$.
\item[(2)] The attractor ${\frak M}_\Lambda$ is invariant in the sense that ${\frak S}({\frak M}_\Lambda) = {\frak M}_\Lambda$. In fact, ${\frak M}_\Lambda$ is the maximal invariant compact set in ${\frak X}_\Lambda$. Also, ${\frak M}_\Lambda$ is an invariant compact subset
of the global attractor ${\frak M}$ of the unrestricted dynamics $({\frak X}, \hbox{\sl dist}, {\frak S})$.
\item[(3)] Through every point $(x(0), w)\in {\frak M}_\Lambda$ passes a complete trajectory, i.e., there exist a two-sided sequence of points $\dots, x(-1), x(0), x(1), \dots$ and a two-sided symbolic sequence
$\dots w(-1)w(0)w(1)\dots$ extending $w$ (in the extension of the subshift $\Lambda$) such that
$S_{w(i)}(x(i)) = x(i+1)$ for all integers $i$.
\item[(4)] Let $K_\Lambda$ denote the projection of the attractor ${\frak M}_\Lambda$ onto the $X$ component. The set $K_\Lambda$ is a compact subset
of the set $K$ of Theorem \ref{two}.
There exist compact sets $A_0,\dots , A_{N-1}$ such that
 $K_\Lambda = A_0\cup A_1\cup\dots\cup A_{N-1}$ and
\begin{equation}\label{frac-at}
K_\Lambda = A_0\cup A_1\cup\dots\cup A_{N-1} =
S_0(A_0)\cup S_1(A_1)\cup\dots\cup S_{N-1}(A_{N-1})
\end{equation}
\item[(5)] In general, the attractor ${\frak M}_\Lambda$
is not a product. There may be infinitely many different sets among the slices ${\cal M}(w) = \{x\in X\,:\,(x, w)\in {\frak M}_\Lambda\}$. However, if the subshift $\Lambda$ is sofic, the number of different slices is finite.
\end{enumerate}
\end{theorem}
\noindent This theorem is proved in section \ref{rd}. For information on sofic shifts see \cite{Lind-Marcus} and section \ref{rd}.

Restricted dynamics of a sort has been considered previously, see \cite{MW,MU}. For example, the graph directed Markov systems of \cite{MU} describe
iterations of uniformly contracting maps indexed by the edges of a directed (possibly infinite) graph.
In this case there is a correspondence between the
points of the limit set and the infinite walks
through the graph (the coding space).
Similarly, the directed IFSs discussed in \cite{Barnsley-2} are defined with the help of the
aforementioned correspondence, and the fractal (or attractor) $K_\Lambda$ is understood in terms of the map from the code space to $K$ as the image of $\Lambda$,
\cite[Theorem 4.16.3]{Barnsley-2}. The correspondence between the points of the code space, $\Sigma$, and the points of $K$ is possible because
the maps are contractions (right away, or eventually). Our approach gives a new
and more general
view on restricted dynamics. We
justify the name -- attractor -- and unveil attractors'
more subtle structure
(assertion 5). This new approach allows us to
work in a much more general setting and with transformations that are not contractions.
We do not have and do not use a map from the code
space into the attractor.

\medskip

We should mention the paper of Andres and Fi\v ser,
\cite{AF-1}. They use their result of \cite{AF-2}
on the existence of the fractal (the set $K$ in our notation) for an IFS with compact operators $S_j$
to conclude that fixed time solution operators of systems of ordinary differential equations could play the role of maps generating the IFSs. As an illustration they
use five two-dimensional systems of ODEs to produce five
operators (incidentally, contractions, as noted in \cite{AF-1}) and plot the corresponding dragon-tail-like
fractal set. Although their message is that IFSs and fractals can be generated by solution operators of ODEs, their examples can serve as an illustration
for our dynamics with choice attractors (due to
Theorem 2(3)).

\begin{figure}[h]
\centering
\includegraphics[width=2in]{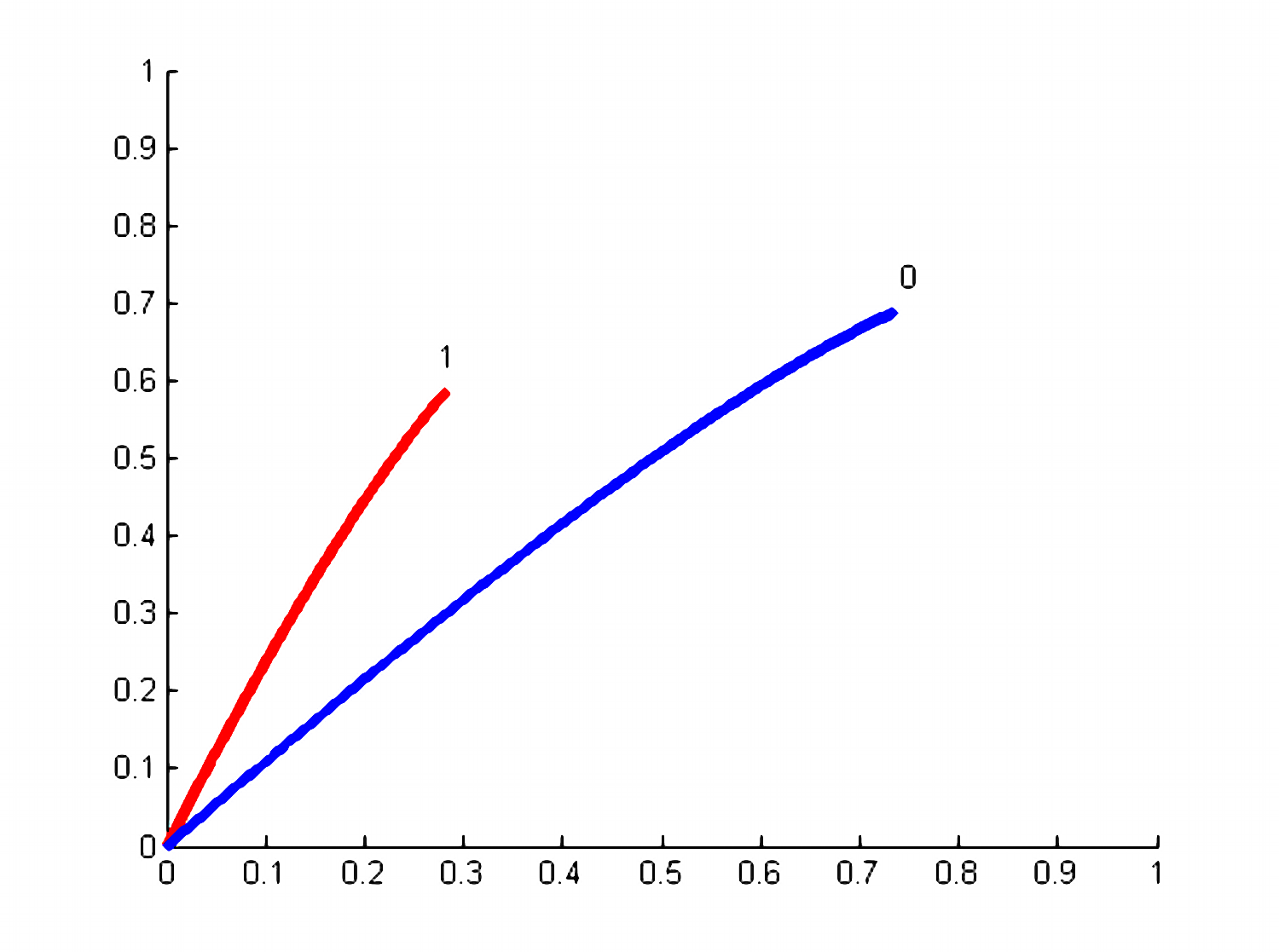}
\caption{Attractors for $(X, d,\,S_0)$ (right) and $(X, d,\,S_1)$ (left).}
\label{f2}
\end{figure}
In section \ref{ex} we apply the theory to a specific example of a discrete Ross-Macdonald type
model of malaria transmission. The model can be viewed as a time discretization (with time step
$\Delta t$) of the ODE model, or as a pre-ODE form
of the model.
The reason we have chosen
this model is because it is simple and we
can visualize all the attractors. The state space is
the unit square $X = \{(x, y)\,:\,0\le x\le 1,\,0\le y\le 1\}$.
We use two sets of parameters, which define two
operators, $S_0$ and $S_1$. Those operators are not contractions, but they are compact, because the system
is finite-dimensional. The discrete dynamical system
generated on $X$ by  $S_0$ has two fixed points,
$(0,0)$, which is unstable, and $(11/15,\,11/16)$,  which is stable. The system generated
by $S_1$ also has two fixed points, $(0,0)$ (unstable) and
$(7/25,\,7/12)$ (stable). The
attractors of the systems
$(X, d,\,S_0)$ and $(X, d,\,S_1)$ are just the heteroclinic trajectories connecting the unstable
and stable fixed points. They are depicted
on figure \ref{f2}.
The effects of  freedom of choice on the dynamics are as follows.
\begin{figure}[h]
\centering
  \begin{minipage}[t]{0.4\linewidth}
     \includegraphics[width=2in]{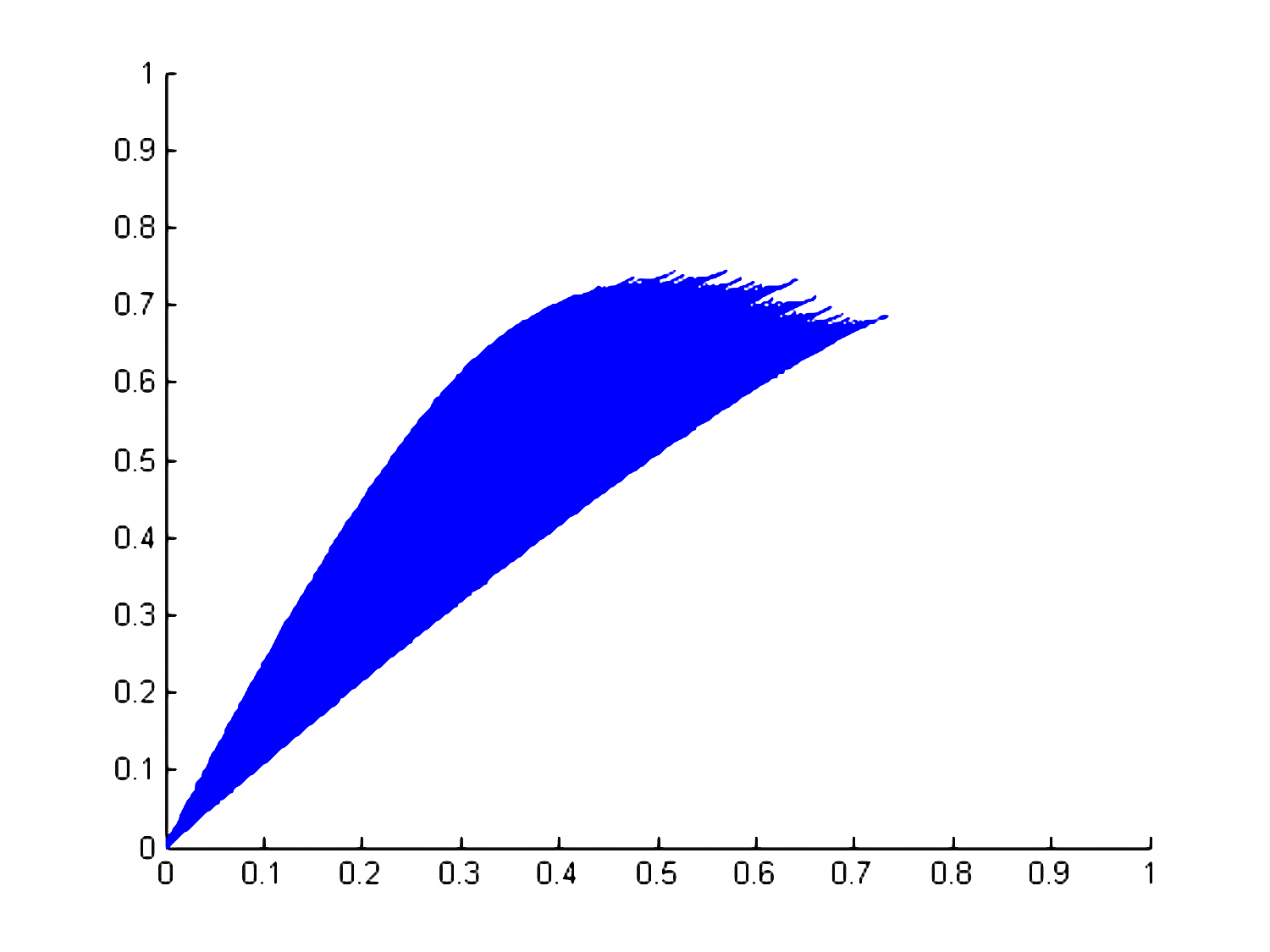}
     \caption{The attractor slice $K$; $\Delta t = 0.05$.}\label{f3}
  \end{minipage}%
  \hspace{0.5in}%
  \begin{minipage}[t]{0.4\linewidth}
     \includegraphics[width=2in]{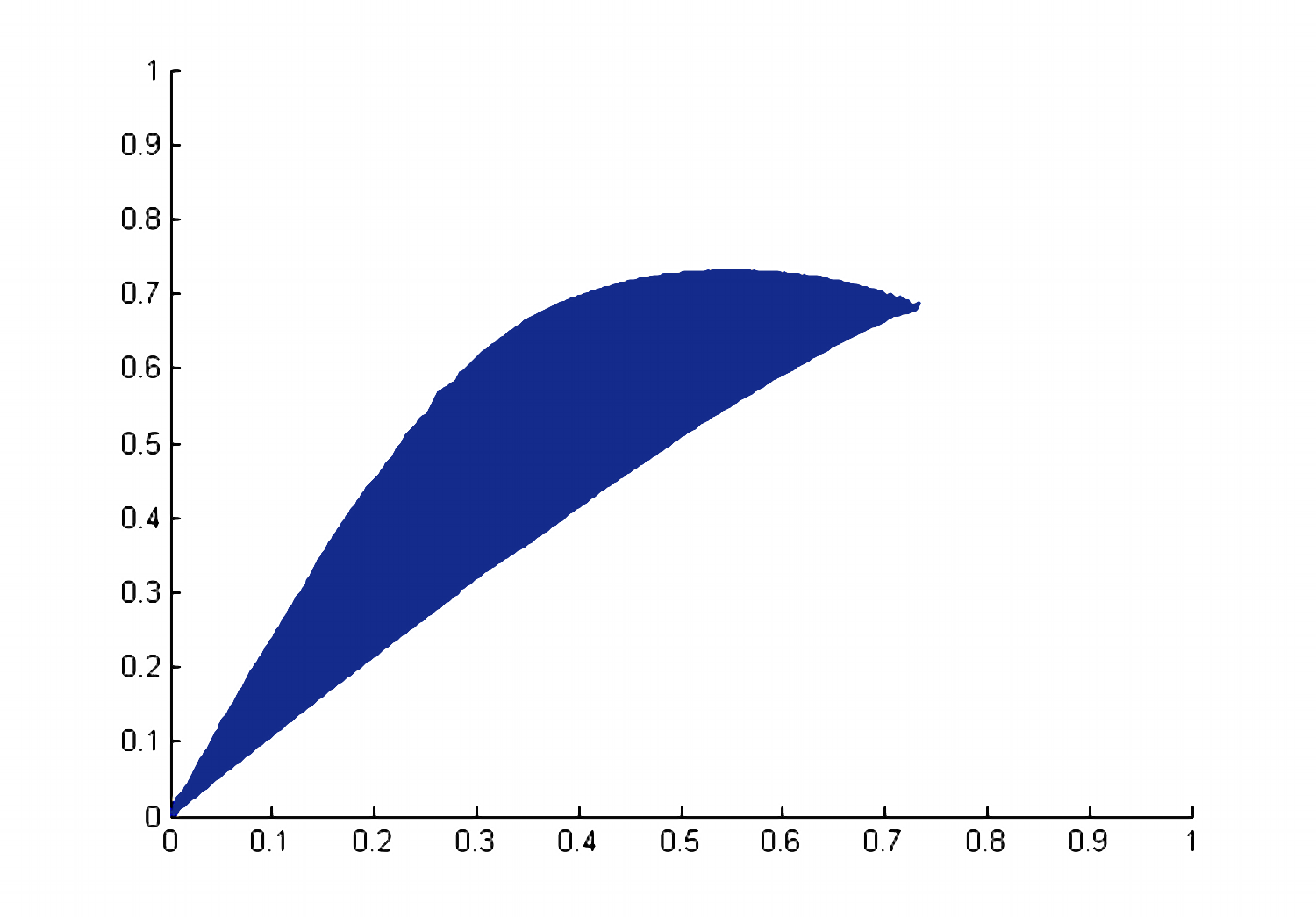}
     \caption{The attractor slice $K$; $\Delta t = 0.005$.}\label{f4}
  \end{minipage}
\end{figure}
For the unrestricted dynamics with choice,
we see (figures \ref{f4} and \ref{f5})
that the attractor, $K$, is a rather big set
with the two individual attractors corresponding to
$S_0$ and $S_1$, respectively, forming parts of the
boundary of $K$ (the right and left sides).
The remaining part of the boundary
is quite irregular when $\Delta t$ is relatively
large, figure \ref{f3}. For smaller $\Delta t$, this part of the boundary becomes much smoother and looks like a smooth curve,
figure \ref{f4}. In the limit $\Delta t\to 0$, the set
$K$ retains its two-dimensional fullness. It is not
an attractor of an ODE with averaged parameters.
In fact, all such attractors are one-dimensional
(each is a heteroclinic trajectory connecting two fixed points) and lie inside of $K$.
\medskip

We consider also the  dynamics with restricted
choice corresponding to the golden mean subshift.
We exhibit two different slices in the global
attractor. Their projections on the state space
do overlap, and their union is smaller than the set
$K$ for the full shift.
\medskip

Finally, we remark that more general compact and condensing
operators will be needed in the study of
dynamics with choice
related to nonlinear dissipative partial differential equations, which we plan to address
at a later time.

\medskip

{\bf The structure of the paper.} This Introduction is followed by two chapters. In the long
chapter \ref{theory} we address theoretical questions. In section \ref{attr} we give
definitions and state the basic results
related to global compact attractors. Section
\ref{mnc} deals with measures of noncompactness.
In section \ref{dwc} we study dynamics with choice.
Dynamics with restricted choice is studied
in section \ref{rd}. In chapter \ref{ex}
we analyze a simple example illustrating some of
our theoretical results. Despite its simplicity,
this example shows that the result of dynamics with choice is larger than the sum of its parts.


\section{Theory}\label{theory}

\subsection{Attractors: general facts}\label{attr}
We start by collecting the basic facts about attractors. There are several books
such as \cite{BV,Hale,Lad-1,Lad-2,SY} devoted to this subject. Our presentation is closer to
\cite{Lad-1}. We present only the results that we need. For the proofs of Theorems \ref{gca} and
\ref{prop} see the books quoted above.

Let $Y$ be a complete metric space with metric
$\d$ and let
$\Phi:\,Y\to Y$ be a continuous map. Iterations,
$\Phi^n$, of
$\Phi$ define a discrete (semi)dynamical system
on $(Y, \d )$. It is useful to consider not only
the dynamics of individual points under
the action of $\Phi$, but, more generally, the dynamics
of bounded sets. Denote by ${\cal B}(Y)$ the collection of all bounded subsets of $Y$.
We say that the set $A\in {\cal B}(Y)$ attracts
the set $B\in {\cal B}(Y)$ if
\[
\vec\d \,(\Phi^n(B),\,A)\,\underset{n\to\infty}\to 0\,,
\]
where the one-sided distance between two sets, $\vec\d \,(C,\,A)$, is understood as
$\sup\limits_{y\in C}\,\d\,(y,\,A)$.

\begin{defin}\label{glob-at}
We call a set $\frak M\subset Y$ the global compact attractor of the system $(Y, \d , \Phi)$ if
\begin{itemize}
\item $\frak M$ is compact,
\item $\frak M$ attracts every bounded subset of $Y$,
\item $\frak M$ is the minimal set with these two properties.
\end{itemize}
\end{defin}
\noindent

\noindent For a system to possess a global compact attractor,
it should enjoy certain properties, namely, some form of  compactness and some dissipativity.
Here is the basic existence (and uniqueness) result.

\begin{theorem}\label{gca}
The semidynamical system
$(Y, \d , \Phi)$ has a global compact attractor
if and only if it enjoys the following two properties:
\begin{enumerate}
\item (``compactness") For every bounded sequence
$(y_k)$ in $Y$ and every increasing sequence
of integers $n_k\to +\infty$, the sequence
$\Phi^{n_k}(y_k)$ has a convergent subsequence.
\item (``dissipativity") There exists a bounded set
${\bf B}\subset Y$ which absorbs
every bounded set in the sense that for every $A\in{\cal B}(Y)$
there exists $m(A)>0$ such that $\Phi^n(A)\subset {\bf B}$ for all $n\ge m(A)$.
\end{enumerate}
\end{theorem}

\noindent Some of basic properties of a global compact attractor are collected in the following theorem.

\begin{theorem}\label{prop}
Assume that $\frak M$
is the global compact attractor of the semidynamical system $(Y, \d , \Phi)$. Then
\begin{enumerate}
\item $\frak M$ is the union of all possible
limits of sequences of the form $\Phi^{n_k}(y_k)$,
where $y_k$ is a bounded sequence in $Y$ and $n_k\to \infty$.
\item $\frak M$ is (strictly) invariant:
$\Phi ({\frak M}) = {\frak M}$.
\item $\frak M$ is the union of all closed bounded sets $A$ with the property $A\subset \Phi(A)$.
\item $\frak M$ is the maximal closed set
with the property $A\subset \Phi(A)$; in particular,
$\frak M$ is the maximal (strictly) invariant closed  set.
\item Through every point $y\in {\frak M}$
passes a complete trajectory, i.e., there exists
a two-sided sequence $\dots, y_{-2}, y_{-1}, y_0, y_1,\dots $ of points in $\frak M$ such that $y_0 = y$ and
$ y_{m+1} = \Phi (y_m)$ for all integers $m$.
\item $\frak M$ is the union of all complete, bounded trajectories in $Y$.
\end{enumerate}
\end{theorem}

\noindent In applications, people do not  verify  the ``compactness" property of Theorem \ref{gca} directly. Instead, they use one of the known sufficient conditions that imply it. Two of the most
useful sufficient conditions are:
\begin{itemize}
\item $\Phi$ is a compact map (i.e., $\Phi: Y\to Y$ is continuous and maps bounded sets into relatively compact sets);
\item in the case $Y$ is a Banach space,
$\Phi$ is a sum of a compact operator and a strict contraction.
\end{itemize}
Compact $\Phi$ arise, e.g., in the finite-dimensional dynamics described by differential
or difference equations, or, in the infinite dimensional case, in dynamics described by
parabolic equations. The ``compact + contraction"
$\Phi$ appear, e.g., in hyperbolic problems with damping. Each of the two sufficient conditions implies that $\Phi$ is condensing with respect
to some measure(s) of noncompactness. Since
measures of non-compactness and condensing
operators are not widely known, below we give
a brief account of the facts we need
and refer to \cite{AKPRS} for more details.


\subsection{Measures of noncompactness}\label{mnc}

Measures of noncompactness assign real non-negative numbers to bounded sets with value $0$ assigned
exclusively to relatively compact sets. The basic
examples are the Kuratowski measure of noncompactness $\alpha$ and the Hausdorff
measure of noncompactness $\chi$. By definition,
$\alpha(A)$ is the infimum of numbers $\epsilon > 0$
such that $A$ admits a finite cover by sets of diameter less than $\epsilon$. The number
$\chi(A)$ is the infimum of those $\epsilon > 0$
for which $A$ possesses a finite $\epsilon$-net in $Y$. In this paper we adopt the following definition
of a general measure of noncompactness (our definition  differs from that in \cite{AKPRS}).
\begin{defin}\label{nc}
A function $\psi$ assigning non-negative real numbers to bounded subsets of (a complete metric)
space $Y$ will be called a measure of noncompactness iff it has the following properties:
\begin{itemize}
\item[(i)] $\psi(A) = 0$ if and only if $A$ is relatively compact;
\item[(ii)] If $A_1\subset A_2$, then $\psi(A_1)\le \psi(A_2)$\ ;
\item[(iii)]  $\psi(A_1\cup A_2) = \max\;\{\psi(A_1),\,\psi(A_2)\}$\ ;
\item[(iv)] There exists a constant $c(\psi)\ge 0$
such that
\[
|\psi(A_1) - \psi(A_2)|\,\le\,c(\psi)\,d_H(A_1, A_2)\,,
\]
where $d_H$ is the Hausdorff distance,
\[
d_H(A_1, A_2) =
\max \;\{ \vec{\hbox{\sl dist}}\,(A_1, A_2),\,\vec{\hbox{\sl dist}}\,(A_2, A_1)\}\,.
\]
\end{itemize}
\end{defin}
\noindent Both $\alpha$ and $\chi$ enjoy all these properties.
Note that property (iv) implies that the measures of
noncompactness of a bounded set and its closure are equal:
\begin{itemize}
\item[(v)] \ \ $\psi (\overline{A}) = \psi(A)$\ .
\end{itemize}
\begin{defin}\label{cond}
A continuous bounded map $\Phi :\,Y\to Y$ is called condensing with respect to the measure of noncompactness $\psi$
(we also say $\Phi$ is $\psi$-condensing) iff
$\psi(\Phi(A)) \le \psi(A)$ for any bounded $A$,
and
$\psi(\Phi(A)) < \psi(A)$ if $\psi(A) > 0$
(i.e., if $\overline{A}$ is not
compact).
\end{defin}


\begin{theorem}\label{cond-0}
Consider the system $(Y, \hbox{\sl dist},\,\Phi)$.
Assume that $\Phi$ is condensing with respect to some measure of noncompactness $\psi$ and that
there exists a bounded set ${\bf B}$ which absorbs every bounded set.
Then $(Y, \hbox{\sl dist},\,\Phi)$ possesses a
global compact attractor.
\end{theorem}
In the case $\psi$ is the Kuratowski measure of noncompactness this result is proved by
\v Seda, \cite{Seda}. In Lemma \ref{key} below
we establish a more general result.


\subsection{Dynamics with choice}\label{dwc}

\subsubsection{Words, strings}\label{sigma}
Fix an integer $N > 1$. Using the integers $0,\dots, N-1$ as the {\it alphabet},  construct strings
(words) of finite length and (one-sided) strings of infinite length.
Denote by $\Sigma^*$ the set of all finite length strings (words), and denote by $\Sigma$ the set of all (one-sided) infinite strings.
The word of length $0$ is the empty word.
The set of non-empty words is denoted $\Sigma^+$. Given a string $w\in \Sigma^*\cup\Sigma$, $w(0)$ is  the {\it first} letter of $w$, and $w(k)$ is
the $(k+1)$-st letter of $w$.
The length of $w$ is denoted $|w|$.
If $w$ is a finite string and $u\in \Sigma^*\cup\Sigma$, their concatenation is denoted $w.u$;
if $|w| = m$, then $(w.u)(m + k) = u(k)$ for $k = 0, 1, \dots$. For a $w\in\Sigma^*$ and
$s\in \Sigma^*\cup\Sigma$, we write
$w\sqsubset s$ if $w$ is the beginning of the string $s$, i.e., if there exists  $u\in \Sigma^*\cup\Sigma$ such that $s = w.u$. For an infinite string, $s$,
its first $n$ letters form a word denoted $s[n]$, i.e., $s[n]=s(0)s(1)\dots s(n-1)$. The set of all words of length $m$ will be denoted by $\Sigma^*_m$.

Equip the space $\Sigma^*\cup\Sigma$ with the metric
$d_\Sigma$, where
$d_\Sigma(u,\,v) = 2^{-m}$
if $u[m-1] = v[m-1]$ but  $u[m] \neq v[m]$. It is well-known, \cite{Barnsley-2}, that both $\Sigma^*\cup\Sigma$ and
$\Sigma$ with metric $d_\Sigma$ are compact.
The shift operator, $\sigma$, acts on infinite
strings by deleting the first letter, i.e.,
$\sigma(u) = u(1)u(2)\dots$. The shift operator
maps $\Sigma$ onto itself. It is continuous; in fact,
$d_\Sigma (\sigma(u),\,\sigma(v)) \le 2\,
d_\Sigma (u,\,v)$.


\subsubsection{The skew-product dynamics}\label{skew}

Let $X$ be a complete metric space with metric $d$,
and let
$S_0, S_1,\dots, S_{N-1}$ be continuous, bounded maps
$X\to X$. Define the product metric space
${\frak X} = X\times\Sigma$ with metric $\hbox{\sl dist}$,
\[
\hbox{\sl dist}\,\left((x, u), (y, v)\right) =
d\,(x, y) + d_\Sigma(u, v)\,.
\]
The skew-product dynamics on ${\frak X}$ is generated by the map ${\frak S}:\,{\frak X}\to {\frak X}$ acting according to the rule
\[
{\frak S}\,(x, u) = (S_{u(0)}\,(x),\,\sigma(u))\,.
\]
This map is obviously continuous and bounded.
Because we will consider iterations of ${\frak S}$
such as
\[
{\frak S}^n\,(x, u) =
(S_{u(n-1)}\circ\cdots \circ S_{u(1)}\circ S_{u(0)}\;(x),\,\sigma^n(u))\,,
\]
we introduce the notation
\[
S_w = S_{w(n-1)}\circ\cdots \circ S_{w(1)}\circ
S_{w(0)}\,,
\]
where $w$ is a word of length $n$. Thus, we can write
${\frak S}^n\,(x, u) = (S_{u[n]}\,(x), \sigma^n(u))$.

\begin{ass1}
Assume there is a closed, bounded set ${\bf B}\subset X$ such that for every bounded $A\subset X$ there exists $m(A)>0$ such that
$S_{w}\,(A) \,\subset {\bf B}$
for every word $w$ of length $n\ge m(A)$.
\end{ass1}
[In applications ${\bf B}$ is usually a closed ball of radius that depends on the parameters of the model. Showing that for different values of
the parameters there is a common  estimate on the radius is enough to verify Assumption 1.]

Let $\psi$ be a measure of noncompactness as in
Definition \ref{nc}.

\begin{ass2}
Assume that each operator
$S_j$ is $\psi$-condensing.
\end{ass2}

\bigskip

We are going to apply Theorem \ref{gca} to prove
the existence of the attractor.  For this to work we need
to justify the following fact:

{\it For every bounded sequence $(x_k, u_k)\in {\frak X}$
and every increasing sequence of integers
$n_k\to +\infty$, the sequence
${\frak S}^{n_k}(x_k, u_k)$ has a convergent
subsequence.}

Thus, pick a bounded sequence $x_k\in X$ and any
sequence $u_k\in\Sigma$. Pick an increasing to
$+\infty$ sequence $n_k$. Since $\Sigma$ is
compact, the sequence $\sigma^{n_k}(u_k)$ has a
convergent subsequence. So, we may assume from the very beginning that
$\sigma^{n_k}(u_k)\to s\in\Sigma$. Denote
$w_k = u_k[n_k]$. We have
\[
{\frak S}^{n_k}(x_k, u_k) = (S_{w_k}\,(x_k),\,\sigma^{n_k}(u_k)).
\]
Since the $\Sigma$-component converges, we need to show that the sequence $S_{w_k}\,(x_k)$ in $X$ has
a convergent subsequence. Because $\Sigma^*\cup\Sigma$ is compact, we can choose a convergent subsequence from $w_k$. We will assume that
$w_k$ itself converges to some $w\in\Sigma$.

\begin{lemma}\label{key}
Under the Assumptions 1 and 2, let $w_n$ be a sequence of finite words of increasing lengths
$|w_n|\to +\infty$ and $w_n\to w\in\Sigma$.
For any bounded sequence $x_n$ the sequence
$S_{w_n}(x_n)$ has a convergent subsequence.
\end{lemma}

\Proof Pick a bounded sequence $x_k$.
By Assumption 1, when the length of the word $w_n$ is sufficiently large, $S_{w_n}(\{x_k\})\subset {\bf B}$. Dropping the first few terms if necessary, we will assume that
$S_{w_n}(\{x_k\})\subset {\bf B}$ for all $n$.
Next, since $w_n\to w$, for arbitrarily large $m_*$
we have $w_n[m_*] = w[m_*]$ for all sufficiently large $n$. Choosing $m_*$ large enough, we have
$S_{w[m_*]}(\{x_k\})\subset {\bf B}$. Writing
$w_n = w[m_*].s_n$, we see that
\[
S_{w_n}(\{x_k\}) = S_{s_n}\left(S_{w[m_*]}(\{x_k\})\right)
\]
and $s_n\to s = \sigma^{m_*}(w)$. Use
Assumption 1 to find $m({\bf B})$ and assemble the
set
\[
\tilde{\bf B} = \bigcup\limits_{v\in \Sigma^*_{m({\bf B})}}\;S_v ({\bf B})\,.
\]
The set $\tilde{\bf B}$
has the property that $S_v(\tilde{\bf B})\subset
{\bf B}$ for all finite words $v$.
If we choose $m_*$ above sufficiently large
(first choose it to guarantee $S_{w[m_*]}(\{x_k\})\subset {\bf B}$ and then increase it by
$m({\bf B})$),
the set $S_{w[m_*]}(\{x_k\})$ will be inside
of $\tilde{\bf B}$. We will reformulate our problem now. We have a sequence $y_k = S_{w[m_*]}(x_k)$
in $\tilde{\bf B}$ and a sequence $s_k\to s$.
We need to show that the sequence
$S_{s_k} (y_k)$ is relatively compact.

 Consider the
positive (infinite) trajectory of $\tilde{\bf B}$ under
the Hutchinson-Barnsley evolution:
\[
C_0 = \tilde{\bf B} \cup \bigcup\limits_{n\ge 1}
\bigcup\limits_{v\in \Sigma^*_{n}}\,
S_v (\tilde{\bf B})\,.
\]
Define inductively
\[
C_{n+1} = \bigcup\limits_{j = 0, \dots, N-1}\,
S_j\,(C_n)\,.
\]
We have $C_0\supset C_1\supset C_2\supset\dots$.

Introduce a collection $\frak H$ of all sets
$A\subset {\bf B}$ that can be represented in the form
\[
A = \bigcup\limits_{n\ge 0} A_n\,,\quad
\textrm{where $A_n$ is a {\it finite} (or empty)  subset of $C_n$}\,.
\]
We show next that every set $A\in {\frak H}$ is relatively compact. Because the sequence
$\{S_{u_n}(y_n)\}$ is in ${\frak H}$, our lemma will be proved. The argument that follows is a modification of a part of \cite[Lemma 1.6.11]{AKPRS}.

First, we claim that there exists a set $A^*\in{\frak H}$ such that
\[
\psi (A^*) = \sup\limits_{A\in{\frak H}}\,\psi(A)\,.
\]
This follows from Lemma 1.6.10 of  \cite{AKPRS}.
Although their lemma is  stated for the Hausdorff
measure of noncompactness,
their proof uses only the properties
(i), (ii), and (iii) of Definition \ref{nc} and,
therefore, works for our $\psi$.

Now, $A^* = \bigcup\limits_{n\ge 0} A^*_n$, where
each $A^*_n$ is a finite (or empty) subset of $C_n$.
For every $p\in A^*_n$ with $n\ge 1$ pick
a $q\in C_{n-1}$ and a $j\in [0, 1,\dots , N-1]$
such that $S_j(q) = p$. Denote the resulting subset
of $C_{n-1}$ by
$A_{n-1}$ and define $A = \bigcup\limits_{n\ge 0} A_n$. Clearly, $A\in {\frak H}$. Hence,
$\psi(A)\le \psi(A^*)$. Now, consider the set
${ F}(A) = S_0(A)\cup\dots\cup S_{N-1}(A)$.
Clearly, ${ F}(A)\in {\frak H}$ and
${ F}(A)\supset \bigcup_{n\ge 1} A^*_n$.
Hence,
\[
\psi (A^*) = \psi (\bigcup_{n\ge 1} A^*_n) \le
\psi ({ F} (A)) = \max \{\psi (S_0(A)),\,\dots,\,\psi (S_{N-1}(A))\}\le \psi (A)\,.
\]
The properties of $\psi$ used here are, from left to right: the first equality uses properties (i) and (iii), the first inequality follows from (ii),
the second equality follows from (iii), the second inequality is due to the assumption that each operator
$S_j$ is $\psi$-condensing.  Now, if $A$ is not relatively compact, we must have
$\psi(S_j(A))<\psi(A)$ for every $j$.
This would imply
$\psi (A^*) < \psi(A)$ which contradicts the fact
that  $\psi (A)\le \psi(A^*)$. Thus,
$\psi (A^*) = \psi(A) = 0$, and so every set in
${\frak H}$ is relatively compact. This concludes
the proof of Lemma.
\bigskip

Applying now Theorems \ref{gca} and \ref{prop},
we immediately obtain Theorem \ref{one}. We next proceed to the proof of Theorem \ref{two}.

\bigskip

\subsubsection{Proof of Theorem \ref{two}}\label{ifs}
Consider the IFS dynamics $(X, d, {\bar F})$.
Having made Assumption 1
we will prove the existence of a global compact attractor $K\subset X$ if we show that
for every bounded sequence $(x_k)\subset X$
and every increasing sequence $n_k\to +\infty$,
the sequence ${\bar F}^{n_k}(x_k)$ is
relatively compact. But this follows from Lemma \ref{key}. Thus, the existence of the global compact attractor for the IFS is established. The global
compact attractor, $K$, comes with all the properties listed in Theorem \ref{prop}.
In particular, $K$ is the maximal compact set invariant under ${\bar F}$.

To prove that ${\frak M} = K\times\Sigma$
we start by showing that the slices of the attractor
corresponding to different strings are all the same,
i.e., the set $\{x\in X\,:\,(x, s)\in{\frak M}\}$
does not depend on $s$.
\medskip

\noindent{\bf All slices are equal}. Recall that every point $(x, s)$ in ${\frak M}$
is a limit of some sequence ${\frak S}^{n_k}(x_k, s_k)$ with bounded $(x_k)\subset X$ and
$\sigma^{n_k}(s_k)$ converging to $s$. As we argued above, we can write
\[
{\frak S}^{n_k}(x_k, s_k) = (S_{w_k}(x_k), \sigma^{n_k}(s_k))\,,
\]
where $w_k$ is a prefix of length
$|w_k| = n_k$ of the string $s_k$,
i.e., $s_k = w_k.\sigma^{n_k}(s_k)$. The sequence
$S_{w_k}(x_k)$ converges to $x$ and
$\sigma^{n_k}(s_k)$ converges to $s$.
The limit of the pair will not change if we replace
$s_k$ by $w_k.s$. Clearly, for any string
$u\in\Sigma$, we have
\[
\lim\;(S_{w_k}(x_k), \sigma^{n_k}(w_k.u)) = (x, u)\,.
\]
This proves that ${\frak M} = A\times \Sigma$,
with a compact set $A\subset X$.

\medskip

Since $\Sigma = 0.\Sigma\cup 1.\Sigma\cup\dots\cup (N-1).\Sigma$ and since ${\frak S}({\frak M}) = {\frak M}$, we get
${\frak S}(A\times \Sigma) =
(S_0(A)\cup S_1(A)\cup\dots\cup S_{N-1}(A))\times \Sigma = A\times \Sigma$. In other words,
$A = S_0(A)\cup S_1(A)\cup\dots\cup S_{N-1}(A)$.
Because $K$ is the maximal compact in $X$ with this property, we have $A\subset K$. On the other hand,
${\frak S}(K\times\Sigma) = K\times\Sigma$. Since
$A\times \Sigma$ is the maximal compact in ${\frak X}$ with this property, we have $K\subset A$,
and hence, $A = K$. This completes the proof of
Theorem \ref{two}.


\subsubsection{Individual attractors}\label{ias}

Every fixed strategy  also generates a dynamics
on $X$: if $w\in\Sigma$ is the (fixed) strategy,
then an $x\in X$ moves to $S_{w(0)}(x)$, then to
$S_{w(1)}\left(S_{w(0)}(x)\right)$, then to
$S_{w(2)}\left(S_{w(1)}\left(S_{w(0)}(x)\right)\right)$, etc. Denote this dynamics by
$(X, d, w)$. This is not a (semi)dynamical system,
but we should not worry about names. Certain
important notions related to the long-term behavior
with natural adjustments still make sense. For example, the individual, i.e.,
corresponding to an individual strategy $w$,
trajectory of a set $B$ is the union
\[
B \cup S_{w[1]}(B) \cup S_{w[2]}(B)\cup\dots\;.
\]
We define the individual $\omega$-limit set of a bounded set $B$ as
\[
\omega(B,\,w) = \{ y\in X\,:\,y = \lim
S_{w[n_k]}(y_k) \quad\textrm{for some sequence
$(y_k)$ in $B$}\}\,.
\]
By analogy with Definition \ref{glob-at}, we say that
a set $A$ is the global compact attractor of
system $(X, d, w)$ if it is
the minimal set with the following two properties:
$A$ is compact and $A$ attracts every bounded set under the strategy $w$,
i.e., for any bounded $B$, we have
 $\lim_{n\to\infty}
 \vec{\hbox{\sl dist}}\,(S_{w[n]}(B),\,A) = 0$.

Next theorem establishes the existence of individual
compact attractors, ${\cal A}_w$, of systems
$(X, d, w)$.
Along the way we establish various properties of
the $\omega$-limiting sets $\omega(B,\,w)$.

\begin{theorem}\label{indiv}
Under the Assumptions 1 and 2, every system
$(X, d, w)$ has the  global compact attractor,
which we denote by
${\cal A}_w$. This attractor is the intersection
of the closures of the tails of the trajectory
of the absorbing set $\tilde{\bf B}$,
\[
{\cal A}_w = \bigcap\limits_{n\ge 1}\;
\bigcup\limits_{k\ge n}
\overline{S_{w[k]}(\tilde{\bf B})}\,.
\]
The attractor, ${\cal A}_w$, is the union of
all $\omega(B, w)$ with bounded $B$.
\end{theorem}

\Proof We use some notation and keep in mind the argument from the proof of
Lemma \ref{key}. Due to Assumption 1, every bounded
set eventually finds itself in the set
$\tilde{\bf B}$ and after that stays there.
 \medskip

\noindent{\bf Step 1. The $\omega$-limit sets of
bounded sets are not empty.}

\noindent Pick a point $x_0\in X$ and follow its trajectory,
$x_n = S_{w[n]}(x_0)$. There will be a time $n$ such
that $x_n\in\tilde{\bf B}\subset C_0$, and then
inevitably
$x_{n+1}\in C_1$, $x_{n+2}\in C_2$, and so on.
By Lemma \ref{key}, the sequence $(x_n)$ is relatively compact. Thus, $\omega(\{x_0\}, w)\neq\emptyset$.
Because $\omega(\{x_0\}, w)\subset \omega(B, w)$
if $x_0\in B$, we have $\omega(B, w)\neq\emptyset$.
\medskip

\noindent{\bf Step 2.  \ $\omega(B, w)$ is the intersection of the closures of the tails of its trajectory, hence $\omega(B, w)$ is closed.}

\noindent Note that $\omega(B, w)$ can be characterized
as follows. $\omega(B,\,w)$ is the set of all
$y\in X$ such that for every $\epsilon>0$
and every integer $k\ge 0$ there exist
and $x\in B$ and $n > k$ so that
$S_{w[n]}(x)\in\;{\cal O}_\epsilon (y)$ (where
${\cal O}_\epsilon (y)$ is the $\epsilon$-neighborhood of $y$).
Yet another way to describe $\omega(B,\,w)$ is
to consider the trajectory of $B$ and its tails:
\[
D_0 = B \cup S_{w[1]}(B) \cup S_{w[2]}(B)\cup\dots\,, \qquad
D_n = \bigcup\limits_{m\ge n} S_{w[m]}(B)\,.
\]
Clearly, $D_0\supset D_1\supset D_2\supset\dots$.
It turns out that
\be\label{omega}
\omega(B,\,w) = \bigcap\limits_{n\ge 0} \,\overline{D_n}\,.
\ee
Indeed, inclusion $\subset$ is obvious. To prove
the ``$\supset$" part, pick a $y$ in the intersection
of the tails and set $\epsilon_n = 2^{-n}$.
In $D_1$ there is a point $y_1 = S_{w[m_1]}(x_1)$, $x_1\in B$,
such that $d\,(y_1,\,y)\le\epsilon_1$.
In $D_{m_1+1}$ there is a point $y_2 = S_{w[m_2]}(x_2)$, $m_2 > m_1$,
such that $d\,(y_2,\,y)\le\epsilon_2$, and so on.
The limit of $y_n$ belongs to $\omega(B,\,w)$, i.e.,
$y\in \omega(B,\,w)$.
\medskip

\noindent{\bf Step 3.  \ $\omega(B, w)$ is compact.}

\noindent Compactness of $\omega(B,\,w)$  will follow
from the fact that the intersection of the
closures of the sets $C_n$ in the proof of Lemma
\ref{key} is compact, because, thanks to Assumption 1, $\bigcap\limits_{n\ge 0} \,\overline{D_n}\,\subset\,\bigcap\limits_{n\ge 0} \,\overline{C_n}$.
Denote $C_* = \bigcap\limits_{n\ge 0} \,\overline{C_n}$.
Since $\omega(B,\,w)$ is not empty,
$C_*$ is not empty as well. And it is closed.
If the set $C_*$ is not compact, then
there exist $\epsilon_0 > 0$ and an infinite sequence $(y_n)\subset C_*$ such that
$d\,(y_n,\,y_m) \ge \epsilon_0$ for all $n$ and $m\neq n$. Since $y_n\in\overline{C_n}$ (in fact, the whole sequence lies in every set
$\overline{C_n}$), there exists a sequence
$y_{nk}\in C_n$ that converges to $y_n$ as $k\to\infty$. For every $\epsilon > 0$ there are numbers
$k_n$ such that $d\,(y_{nk_n},\,y_n) \le \epsilon$
for all $n$. When $\epsilon < \epsilon_0/2$, the Hausdorff distance between the sets
$\{y_{nk_n}\}$ and $\{y_n\}$ is not greater than
$\epsilon$. Using property (iv) of the measure of
noncompactness $\psi$, we obtain
$|\psi\left(\{y_n\}\right) - \psi\left(\{y_{nk_n}\}\right)| \le\; c(\psi)\,\epsilon$. Now,
$\psi\left(\{y_{nk_n}\}\right) = 0$ by Lemma \ref{key}. Then $\psi\left(\{y_n\}\right)\le\; c(\psi)\,\epsilon$. Since this is true for any $\epsilon$,
we obtain
$\psi\left(\{y_n\}\right)= 0$, a contradiction.
This proves that $C_*$ is  compact.

\medskip

\noindent{\bf Step 4.
$\omega(B,\,w)$ attracts $B$.}

\noindent To show that
for every $\epsilon >0$ there exists an $m$ such that $S_{w[n]}(B)\subset {\cal O}_\epsilon
\left(\omega(B,\,w)\right)$ for all $n\ge m$,
we argue by contradiction. Assume there exists
an $\epsilon_0 > 0$ such that
$S_{w[n]}(B)$ does not lie inside
${\cal O}_{\epsilon_0}
\left(\omega(B,\,w)\right)$ for infinitely many $n$.
This means that there is a sequence $x_k$ in $B$
and a sequence $n_k\to +\infty$ such that
$S_{w[n_k]}(x_k)\notin {\cal O}_{\epsilon_0}
\left(\omega(B,\,w)\right)$. But we already know that $S_{w[n_k]}(x_k)$ must have a convergent
subsequence whose limit must be in $\omega(B,\,w)$.
A contradiction.

\medskip

\noindent{\bf Step 5.}\
${\cal A}_w = \omega(\tilde{\bf B},\,w)=
\bigcup\limits_{\textrm{bounded $B$}}\omega(B, w)$.

\noindent Because every bounded set is eventually absorbed by the set $\tilde{\bf B}$, we have
$\omega(B,\,w) \subset \omega(\tilde{\bf B},\,w)$.
Thus, $\omega(\tilde{\bf B},\,w)$ attracts every bounded set.
It is compact and minimal,
hence, it is the global compact attractor. The theorem is proved.

\bigskip


\subsubsection{Interplay between individual attractors}\label{inter}

Recall, that (with Assumptions 1 and 2) the global
attractor $\frak M$ of $(X, \d, \Sigma)$ is a product
${\frak M} = K\times \Sigma$.

We start with a few simple observations.
\begin{lemma}
${\cal A}_w \subset F\left( {\cal A}_{w}\right)\subset K$,
where $F$ is the Hutchinson-Barnsley operator.
\end{lemma}
\Proof  Pick a point, $x$, in ${\cal A}_w$.
Then $x = \lim S_{w[n_k]}(x_k)$ for some
bounded sequence $(x_k)$ in $X$ and $n_k\to \infty$.
Among the last letters of the words
$w[n_k]$ there is at least one, say, $j$, that repeats infinitely many times. Sparse the sequence so that every  $w[n_k]$ has the last letter $j$.
Then,
\[
S_{w[n_k]}(x_k) = S_{j}\left(
S_{w[n_k-1]}(x_k)
\right)\,.
\]
The sequence $S_{w[n_k-1]}(x_k)$
has a convergent subsequence by Lemma \ref{key},
and the limit is in ${\cal A}_w$. Thus,
$x\in S_j\left({\cal A}_w\right)$.
Lemma is proved.
\medskip

\begin{lemma}
${\cal A}_w \subset {\cal A}_{\sigma(w)}$.
\end{lemma}
\Proof  Again, if $x\in{\cal A}_w$, then $x = \lim S_{w[n_k]}(x_k)$. Clearly,
\[
S_{w[n_k]}(x_k) =
S_{\sigma(w)[n_k-1]}\left(S_{w(0)}(x_k)\right)
\,.
\]
The sequence $(S_{w(0)}(x_k))$ is bounded
and $\sigma(w)[n_k-1]\to \sigma(w)$.
Lemma is proved.

\begin{corollary}\label{per-at} If the string $w$ is periodic, then
${\cal A}_w = {\cal A}_{\sigma(w)}$.
\end{corollary}

The union of
individual attractors ${\cal A}_w$ lies inside of $K$,
\be\label{inclusion}
\bigcup\limits_{w\in\Sigma} {\cal A}_w\,\subseteq K\,.
\ee
There are many important cases when this union equals $K$.


\begin{lemma}\label{equality}
We have $\bigcup\limits_{w\in \Sigma} {\cal A}_w = K $  in
each of the following cases:
\begin{enumerate}
\item[a)] Operators $\{S_j\}$ are eventually strict contractions,
i.e., there exist a $0 < \gamma < 1$ and an integer $M \geq 1$
such that for any finite word $w^*$ of length $\ge M$ the operator
$S_{w^*}$ is a contraction with factor $\gamma$. (This condition is automatically satisfied if each $S_j$ is a strict contraction.)
\item[b)] $S_j^{-1}(K)\supseteq K$ for $j = 0,\dots, N-1$.
\item[c)] Each operator $S_j$ is invertible on $K$.
\end{enumerate}
\end{lemma}
\Proof \ The inclusion (\ref{inclusion}) is obvious. To prove the equality in the special cases a) and b),  pick an $x\in K$.
There exists a sequence of points $\{x_k\}\subset K$, and a sequence $w_{n_k}$
of lengths $n_k$ increasing to infinity such that
$x=\lim\limits_{k\rightarrow \infty} S_{w_{n_k}}(x_k)$. We claim that $x \in {\cal A}_u$, where $u=w_{n_1}.w_{n_2}\ldots w_{n_k}\ldots$.
Denote $u[m_k]=w_{n_1}.w_{n_2}\ldots w_{n_k}$. The lengths
of the words $u[m_k]$  go to infinity.

In the case {\it a)},
for every $k$ and any $y\in K$ we have

\[
\ba{ll}
d(S_{w_{n_k}}(x_k),S_{u[{m_k}]}(y))=d(S_{w_{n_k}}(x_k), S_{w_{n_k}}S_{u[{m_{k-1}}]}(y))\\
=d(S_{w_{n_k}}(x_k),S_{w_{n_k}}(z_k))
\ea
\]

\noindent where $z_k=S_{u[{m_{k-1}}]}(y)$.
Then, $d(S_{w_{n_k}}(x_k),S_{w_{n_k}}(z_k))\leq \gamma^{l_k}\;d(x_k,z_k)\leq \gamma^{l_k}\;\hbox{diam}(K)$, where
 $l_k$ is the round down of ${n_k}/{M}$. Therefore, $d(S_{w_{n_k}}(x_k),S_{u[{m_k}]}(y))\rightarrow 0$, as $k\rightarrow \infty$.
Since, $\lim\limits_{k\rightarrow \infty} S_{u[{m_k}]}(y) \in {\cal A}_u$, and
$\lim\limits_{k\rightarrow \infty} S_{u[n_k]}(y)=\lim\limits_{k\rightarrow \infty} S_{w_{n_k}}(x_k)=x$, it
follows that $x\in {\cal A}_u$ and the inclusion $K\subset \bigcup\limits_{w\in \Sigma} {\cal A}_w$ is proved.

In the second case, since $S_j^{-1}(K)\supseteq K$,
for every $y\in K$ there exist
$z_j\in K$ with $y=S_j(z_j)$, $j\in \{0,1,\ldots,N-1\}$. Therefore, for every $k$,
we can find $y_k\in K$ such that $S_{u[{m_{k-1}}]}(y_k)=x_k$.
Then, $S_{u[{m_k}]}(y_k)=S_{w_{n_k}}S_{u[{m_{k-1}}]}(y_k)=S_{w_{n_k}}(x_k)$.
It follows that $x\in {\cal A}_u$.

Finally, {\it c)} is a special case of {\it b)}.
This concludes the proof.

\begin{remark}\label{invertible}
The case {\it c)} may seem too restrictive. However,
there are many situations where the operators $S_j$ are not
invertible on $X$ but are invertible on the attractor $K$.
This was first observed by Ladyzhenskaya in the case of
Navier-Stokes equations, \cite{Lad-0}. The fact is due to
the invariance of $K$ and, what is called, {\it backward uniqueness}
property of certain parabolic-like equations.
\end{remark}
\bigskip

Although $K$ equals the union of individual attractors in many cases,
there are situations when $K$ is strictly larger than that union.
This is what we call a Gestalt effect. This is a new phenomenon.
As we have shown in Lemma \ref{equality}, the Gestalt effect
cannot occur when operators $S_j$ are contractions.

\bigskip

\noindent{\bf Example of a Gestalt effect.}
\medskip

In this example the state space $X$ will be the space
$\Sigma_2$  of one-sided infinite strings of $0$'s and $1$'s.
There will be two operators, $S_0$ and $S_1$, defined as follows:
\[
S_0(v) = v(2).v\,,\quad S_1(v) = v(1).v
\]
for all $v = v(0)v(1)v(2)v(3)\dots\in X$.
The conditions of Theorem \ref{one} are satisfied, so
let ${\frak M} = K\times \Sigma_2$ be the global compact attractor
of the
corresponding dynamics with choice.
[Note that the global compact attractor of the
system generated by
$S_0$ is the set
of all strings with period $3$, and the attractor of the system
generated by $S_1$ is the set
of all strings with period $2$.]

We claim that the sequence
$u=000\overline{100}$ is in $K$ but not in ${\cal A}_w$
for any $w\in\Sigma_2$.
Let $v=001.\sigma^3(v)$, i.e., the first three symbols of $v$ are $001$, and let $w_k=000...0001$
with $3k$ zeros before 1. Then, for every k, $S_{w_k}(v)=0.001001...001$ with $001$
repeating $k$ times. Therefore, $S_{w_k}(v)\rightarrow u$ as $k\rightarrow \infty$, i.e.,
$u\in K$. To show that $u$ does not belong to the union $\bigcup\limits_{w\in \Sigma} {\cal A}_w$,
we argue by contradiction. If $u\in {\cal A}_s$, then
there exists a sequence $v_{k}\in \Sigma_2$  such that
$\lim\limits_{k\rightarrow \infty} S_{s[n_k]}(v_k)=u$,  where $n_k\nearrow \infty$.
Therefore, we can find $l$, such that $S_{s[n_{l}]}(v_{l})$,
$S_{s[n_{l+1}]}(v_{l+1})$, \dots, $S_{s[n_{l+8}]}(v_{l+8})$,
all begin with $000100100...$.
Since
$S_{s[n_{l+1}]}(v_{l+1})=S_{s(n_{l+1})}...S_{s[n_{l}]}(v_{l+1})=0001001001...$, and the action of operators $S_0$ and $S_1$
depends only on the first three symbols in the strings,
it follows that $v_{l}[3]\neq v_{l+1}[3]$, because if
$v_{l}[3] = v_{l+1}[3]$, then $S_{s[n_{l+1}]}(v_{l+1})$
starts with at least 4 zeros, i.e., $0000100100...$, which is impossible. Similarly, $v_{l+k}[3]\neq v_{l+j}[3]$ for $j, k = 0, \dots, 8$, $j\neq k$. But there can be only $8$ different
three-letter words in $2$ symbols. A contradiction.
Hence, $u$ does not
belong to the $\bigcup\limits_{w\in \Sigma} {\cal A}_w$.



\subsection{Dynamics with restricted choice}\label{rd}

As in section \ref{sigma}, $\Sigma$ denotes the space of one-sided infinite
strings on $N$ symbols, and $d_\Sigma$ is the metric on $\Sigma$.
Let $\Lambda$ be a subshift of $\Sigma$, i.e., $\Lambda$ is a closed subset of $\Sigma$ and $\sigma(\Lambda) = \Lambda$. Dynamics with restricted choice is defined on the space ${\frak X}_\Lambda = X\times \Lambda$ by the operator
${\frak S}:\;(x, w)\mapsto (S_{w(0)}(x),\,\sigma(w))$,
where the strings $w$ are now taken from $\Lambda$ only.

We assume that $X$ and $S_0,\dots, S_{N-1}$ satisfy our Assumptions
1 and 2. The existence of the global compact attractor,
${\frak M}_\Lambda$, then follows from the abstract result, Theorem
\ref{gca}. The assertions 2 and 3 of Theorem \ref{three} are
among the general properties of global compact attractors,
see Theorem \ref{prop}. Denote by
$K_\Lambda$ the projection of ${\frak M}_\Lambda$
onto the $X$ component. Clearly,
$K_\Lambda$ is compact. Also, $K_\Lambda$ is a subset of the slice $K$ corresponding to
the full shift $\Sigma$, as in Theorem \ref{two}.
Because of the invariance property of ${\frak M}_\Lambda$, for every point
$y\in K_\Lambda$ there is a $j$, one of the
symbols $0, \dots, N-1$, and a point
$x\in K_\Lambda$ such that $y = S_j(x)$. Define
the sets $A_j = \{x\in K_\Lambda\,:\,S_j(x)\in K_\Lambda\}$. It is easy to see that each
$A_j$ is compact and
$K_\Lambda = A_0\cup A_1\cup\dots \cup A_{N-1}$.
By construction, we have
$A_0\cup A_1\cup\dots \cup A_{N-1} = S_0(A_0)\cup S_1(A_1)\cup\dots \cup S_{N-1}(A_{N-1})$.
\medskip

To analyze the slices ${\cal M}_\Lambda(s) = \{x\in X\,:\,(x, s)\in {\frak M}_\Lambda\}$, we follow the argument of the corresponding
part of section \ref{ifs}.

Every point $(x, s)\in {\frak M}_\Lambda$ is the limit of the form
\[
(x, s) = \lim\limits_{n_k\to \infty}\,\left(S_{w_k}(x_{n_k}), \sigma^{n_k}(s_{n_k})\right)\,,
\]
where $(x_n)$ is a bounded sequence in $X$, $(s_n)$ is a bounded sequence in $\Lambda$, and $w_k$ is the prefix of $s_{n_k}$,
$s_{n_k} = w_k.\sigma^{n_k}(s_{n_k})$. Because ${\frak M}_\Lambda$
is invariant under ${\frak S}$ and we know that the unrestricted
dynamics has the global compact attractor ${\frak M} = K\times \Sigma$,
the sequence $(x_n)$ can be taken from the compact $K$, and we may assume that $x_{n_k}\to x_*\in K$. Also, we may assume that the words $w_k$ converge (to some infinite string $w_*\in\Lambda$). The strings
$\sigma^{n_k}(s_{n_k})$ converge to $s$. Consider all strings
$u\in\Lambda$ such
that $w_k.u$ is a string in $\Lambda$ for infinitely many $k$.
For every such $u$ we will have $x\in M_\Lambda(u)$.

We see that the number of different slices of the attractor
${\frak M}_\Lambda$ may depend on the sequence $x_{n_k}$,
but more importantly, it depends on what strings can be attached to
convergent sequences of finite words in $\Lambda$.
\medskip

With every  sequence $(w_k)$ of finite words in
$\Lambda$ we associate the set
${\frak s}((w_k))$ of one-sided infinite
strings $u\in \Lambda$ such that
$w_{k_\ell}.u\in \Lambda$ for some subsequence
$w_{k_\ell}$. In order to prove the third assertion of Theorem \ref{three} we will show that, if $\Lambda$ is
a sofic shift,
the number of different sets among
all ${\frak s}((w_k))$ is finite. The argument will be similar to the proof of
Theorem 3.2.10 in \cite{Lind-Marcus}.

Recall that $\Lambda$ is a sofic shift if it has a presentation by a finite labeled graph, see \cite{Lind-Marcus}. This means that there is a  directed graph, $G = (V, E)$, with a finite number of vertices, $V$, and  edges,
$E$; the edges are labeled by the symbols $0, 1,\dots, N-1$; from every vertex begins at least one  infinite directed path; the labels of the edges in the infinite directed paths form infinite one-sided
strings that exhaust exactly all strings in $\Lambda$.

\begin{lemma}\label{sofic}
If $\Lambda$ is a one-sided sofic subshift of $\Sigma$, then the number of different sets among
all ${\frak s}((w_k))$ is finite.
\end{lemma}
\Proof\ \  Let $G = (V, E)$ be a labeled graph presenting $\Lambda$. Let
$(w_k)$ be a sequence
of finite words allowed in $\Lambda$. For each word
$w_k$ pick a finite directed path in $G$ presenting it. We can find a subsequence, $(w_{k_\ell})$,
such that all the words $w_{k_\ell}$ have the same
terminal vertex in their presentation. If $T$ is such vertex, then $w_{k_\ell}.u\in\Lambda$ for
all infinite paths $u$ starting at $T$. Because the
number of vertices is finite, we are done.
\hfill$\square$%
\bigskip

\begin{remark} Even if the number of different sets among
all ${\frak s}((w_k))$ is $> 1$, the attractor
${\frak M}_\Lambda$ may be a product,
${\frak M}_\Lambda = K_\Lambda\times \Lambda$,
with the same slice for every string in $\Lambda$.
\end{remark}
Indeed, let $N = 2$ and let $\Lambda$ consist
of the periodic string $u = 100100\dots$ and its shifts $\sigma(u) = 00100\dots$ and
$\sigma^2(u) = 0100\dots$. If $(w_k)$ consists of
words ending in $00$, then the only string that
can be attached to $w_k$ is $u$.
If $(w_k)$ consists of
words ending in $1$, then the only string is
$\sigma(u)$, and for words ending in $10$ the only string is $\sigma^2(u)$. Thus, we have three different sets of the form ${\frak s}((w_k))$.
At the same time, the individual attractors
${\cal A}_u$, ${\cal A}_{\sigma(u)}$, and
${\cal A}_{\sigma^2(u)}$, are all equal, as we
argue in Corollary \ref{per-at}.
\bigskip

One may ask whether ${\frak M}_\Lambda$ is always
a product. The answer is no, as the following
example shows.
\begin{figure}[tp]
\centering
\includegraphics[width=0.47\textwidth]{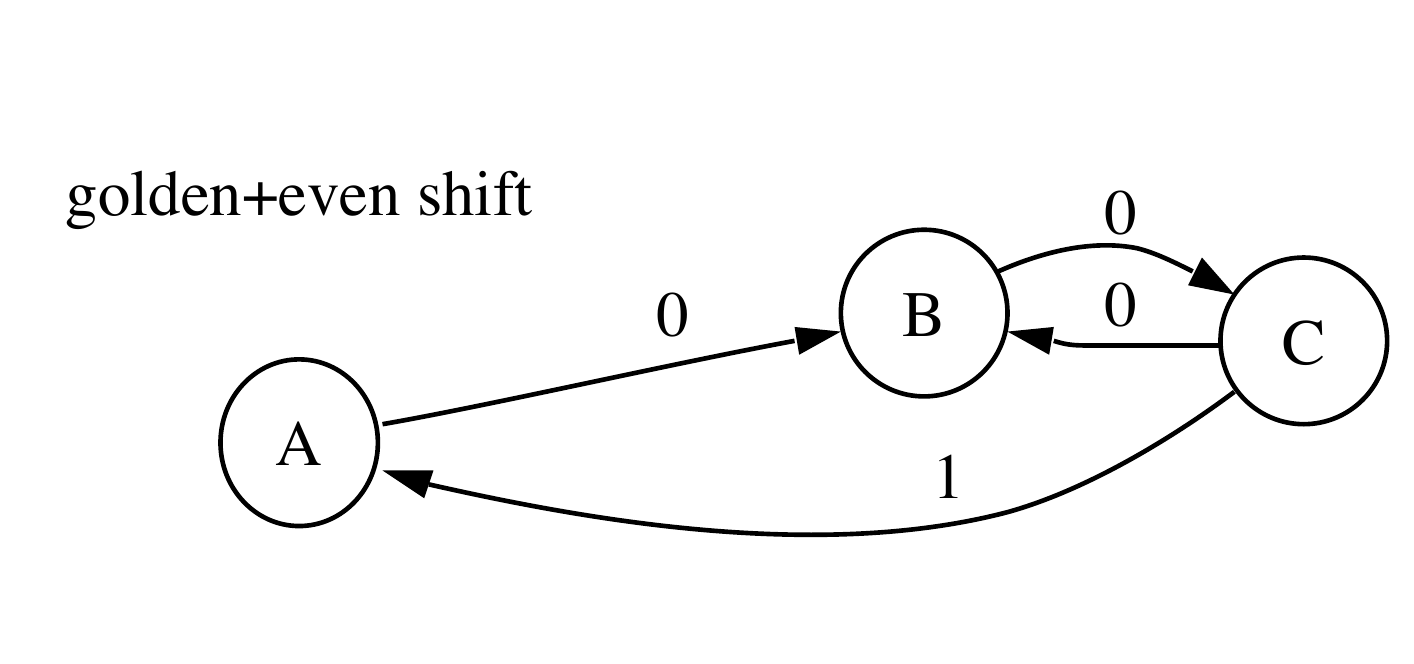}
\hfill
\includegraphics[width=0.47\textwidth]{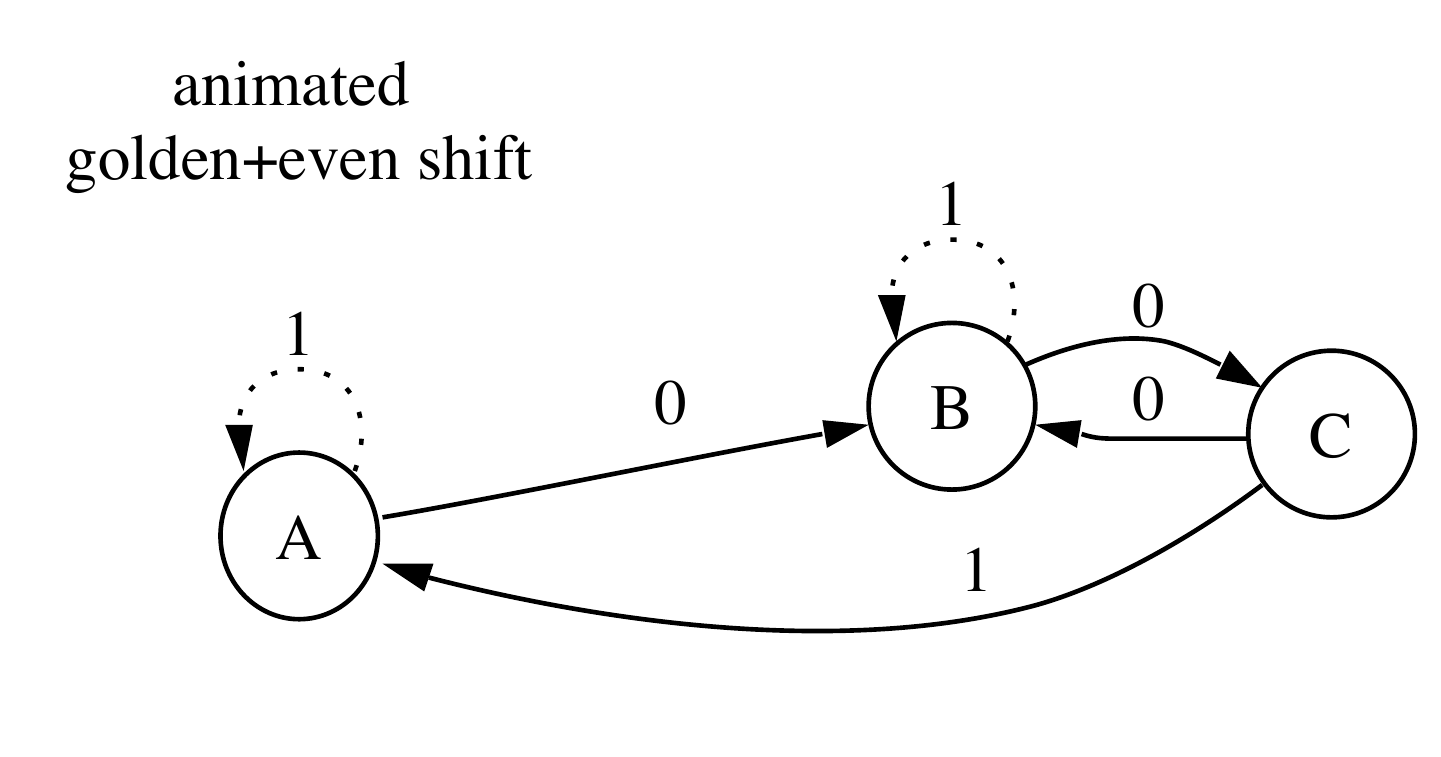}\\
\caption{The golden+even shift and its animation}
\label{f5}
\end{figure}
Let $\Lambda$ be the intersection of the one-sided
golden mean shift with the even shift.
In other words, $\Lambda$ consists
of all sequences of $0$s and $1$s such that between
any two $1$s there are  two or a larger even number
of $0$s. A graph presenting $\Lambda$ is given
on Figure \ref{f5}. We will animate this graph to define the dynamics. First, identify the nodes with three  distinct points $A$, $B$, and $C$
in $\mathbb R^2$, see Figure \ref{f5} left, and
define $X = \{A, B, C\}$. Second,
define the maps $S_0$ and $S_1$ acting on points
as shown by the directed edges labeled correspondingly; for example, $S_0(A) = B$,
$S_1(A) = A$, and $S_0(C) = B$.

Now consider the set $\Lambda^+$ of non-empty
finite words (blocks) of $\Lambda$.
We divide $\Lambda^+$ in three classes and correspondingly divide the strings in $\Lambda$
into three classes. The first class of words in
$\Lambda^+$ consists of the words ending in $1$.
Such words can serve as prefixes of strings starting
with an even (or infinite) number of $0$s.
Denote these classes by $\Lambda^+_A$ and
$_A\Lambda$. The second class of finite words consists of the words ending in odd number of $0$s.
The strings for which such words can serve as prefixes are the strings starting with an odd number
of $0$s. These classes are denoted by
$\Lambda^+_B$ and
$_B\Lambda$. The last class in $\Lambda^+$ consists
of words ending in even number of $0$s. The corresponding strings are those starting with $1$ or
with an even number of $0$s. These are denoted
by $\Lambda^+_C$ and
$_C\Lambda$. By looking at the picture of
the animated shift,
it is easy to identify the possible limits of
sequences
$S_{w_k}(x_k)$ when $w_k$ belong to a particular
class, while $x_k\in \{A, B, C\}$. We see that
if $w_k\in \Lambda^+_A$, then the limit set is
$\{A, B\}$. If $w_k\in \Lambda^+_B$, then the limit set is
$\{B, C\}$. Finally, if $w_k\in \Lambda^+_C$, then the limit set is again
$\{B, C\}$. Thus, there are two different slices
in the attractor ${\frak M}_\Lambda$. One slice is
$\{A, B\}$, and the other is $\{B, C\}$.
We have
${\cal M}_\Lambda(u) =  \{A, B\}$ if $u\in\, _A\Lambda$, and ${\cal M}_\Lambda(u) = \{B, C\}$ if
$u\in\, _B\Lambda\,\cup \,_C\Lambda$.
 The global attractor ${\frak M}_\Lambda$ is a union
of the sets
$\{A, B\}\times _A\Lambda$,  $\{B, C\}\times _B\Lambda$, and $\{B, C\}\times _C\Lambda$.

Another example of different slices appears in
numerical results reported in the next section.



\section{Example}\label{ex}

The simplest mathematical model of malaria transmission goes back to Ross and Macdonald.
The state of the human-mosquito interaction
system is described by the portion of infected
humans, $x$, and the portion of infected
mosquitoes, $y$. The change in time is described
by the following simple system of ordinary differential equations:
\be\label{ode}
\ba{rcl}
\dot x &=& a\,y\,(1 - x) - \;r\,x\\
\dot y &=& b\,x\,(1 - y) - m\,y
\ea
\ee
The nature of the positive
coefficients $a$, $b$, $r$, and $m$ is discussed in \cite{SM}. In particular,
the coefficients $a$ and $b$ are proportional to the biting rate and the transmission efficiencies (infected human to mosquito and  infected mosquito to human),
$r$ is the recovery rate (in humans), and
$1/m$ is the average mosquito life-span. In practice, it is hard to measure these parameters.
Also, there are many factors that affect their values, see \cite{SM}, page 8, and the values may
change in time.

The state space for the model (\ref{ode}) is the
closed square $X = \{(x, y)\,:\,0\le x\le 1,\,
0\le y\le 1\}$. For initial conditions in $X$ the
solution stays in $X$ for all $t$. If the quantity
$R_0 = \frac{ab}{rm}$ is $\le 1$, all trajectories
starting in $X$
converge to the origin, and the global compact attractor consists of a single point, $P_1=(0,0)$.
If $R_0 > 1$, the equilibrium $P_1$ becomes
unstable and there emerges the second fixed point,
$P_2 = (x_*, y_*)$,
inside the square $X$,
\be\label{xstar}
x_* = {ab - rm\over b\,(a + r)}\,,\quad
y_* = {ab - rm\over a\,(b + m)}\,.
\ee
This second equilibrium is stable, and the global
compact attractor of the system consists of the two
equilibria, $P_1$ and $P_2$, and of the heteroclinic
trajectory connecting them (and staying entirely
inside $X$).
The number $R_0$, known as the
basic reproductive number, detects the emergence of
epidemics: when $R_0 > 1$ there is a stable portion of infected population.

We consider a discrete version of equations (\ref{ode}):
\be\label{dode}
\ba{rcl}
x(t + \Delta t) &=& x(t) + \Delta t\,\left(a\,y(t)\,(1 - x(t)) - \;r\,x(t)\right)\\
y(t + \Delta t) &=& y(t) + \Delta t\,\left(b\,x(t)\,(1 - y(t)) - m\,y(t)\right)\,.
\ea
\ee
The time step map $(x(t), y(t))\mapsto (x(t + \Delta t), y(t + \Delta t))$ maps $X$ into itself
provided
\be\label{step}
\Delta t \,<\,\min\,\{\frac{1}{a + r}\,,\;
\frac{1}{b + m}\}\,.
\ee
The fixed points for (\ref{dode}) are the same as for (\ref{ode}). As in the continuous case,
if $ab > rm$ and the time step satisfies (\ref{step}), the global attractor for (\ref{dode})
consists of the two fixed points, $P_1$ and $P_2$,
and the heteroclinic trajectory connecting them.

\begin{figure}[t]
\centering
  \begin{minipage}[t]{0.4\linewidth}
     \includegraphics[width=2in]{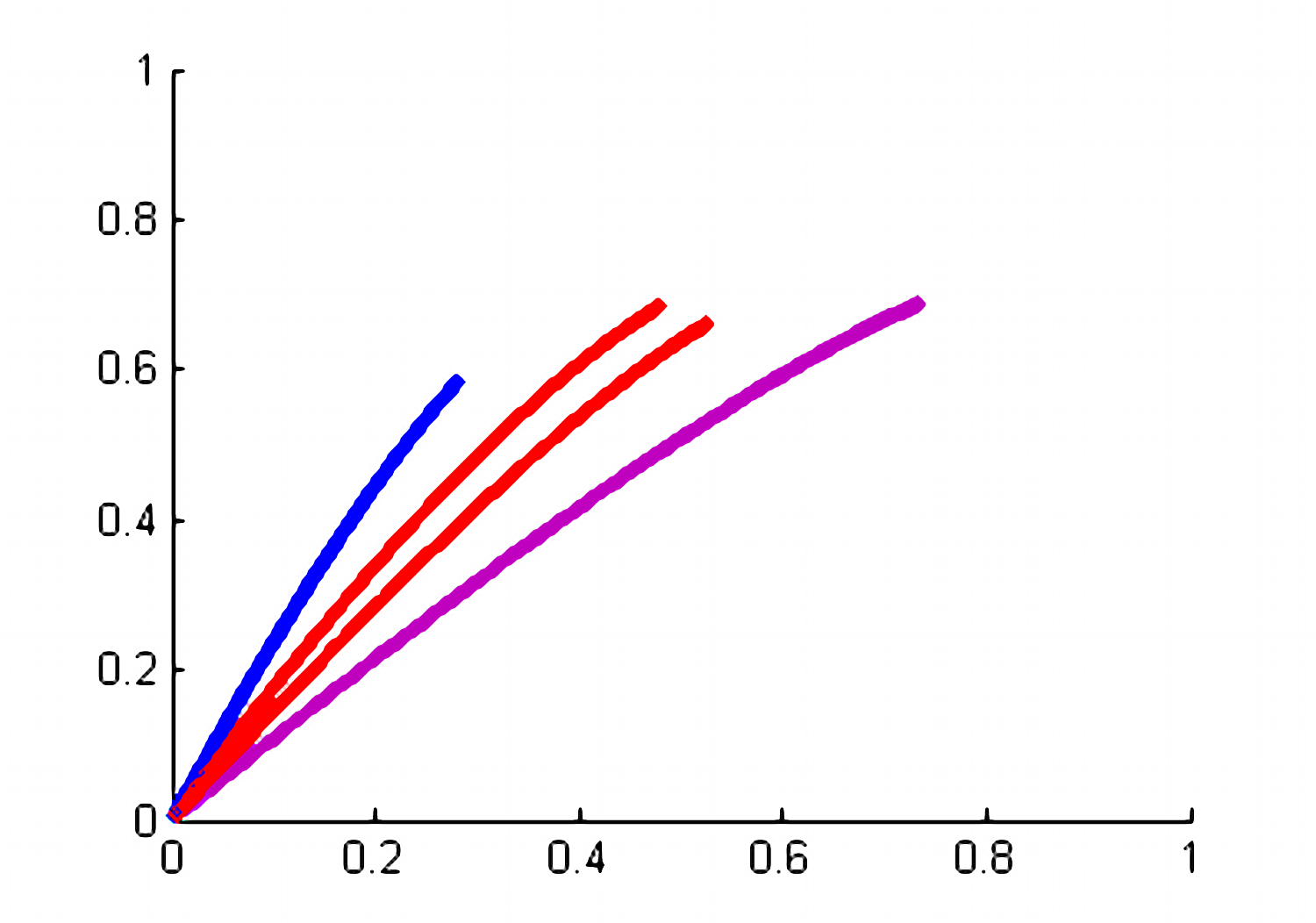}
     \caption{Three individual attractors ${\cal A}_w$ for $\Delta t = 0.05$:  left: $w=111...$; middle two: $w = 1010...$;
     right: $w=000...$.}\label{f6}
  \end{minipage}%
  \hspace{0.5in}%
  \begin{minipage}[t]{0.4\linewidth}
     \includegraphics[width=2in]{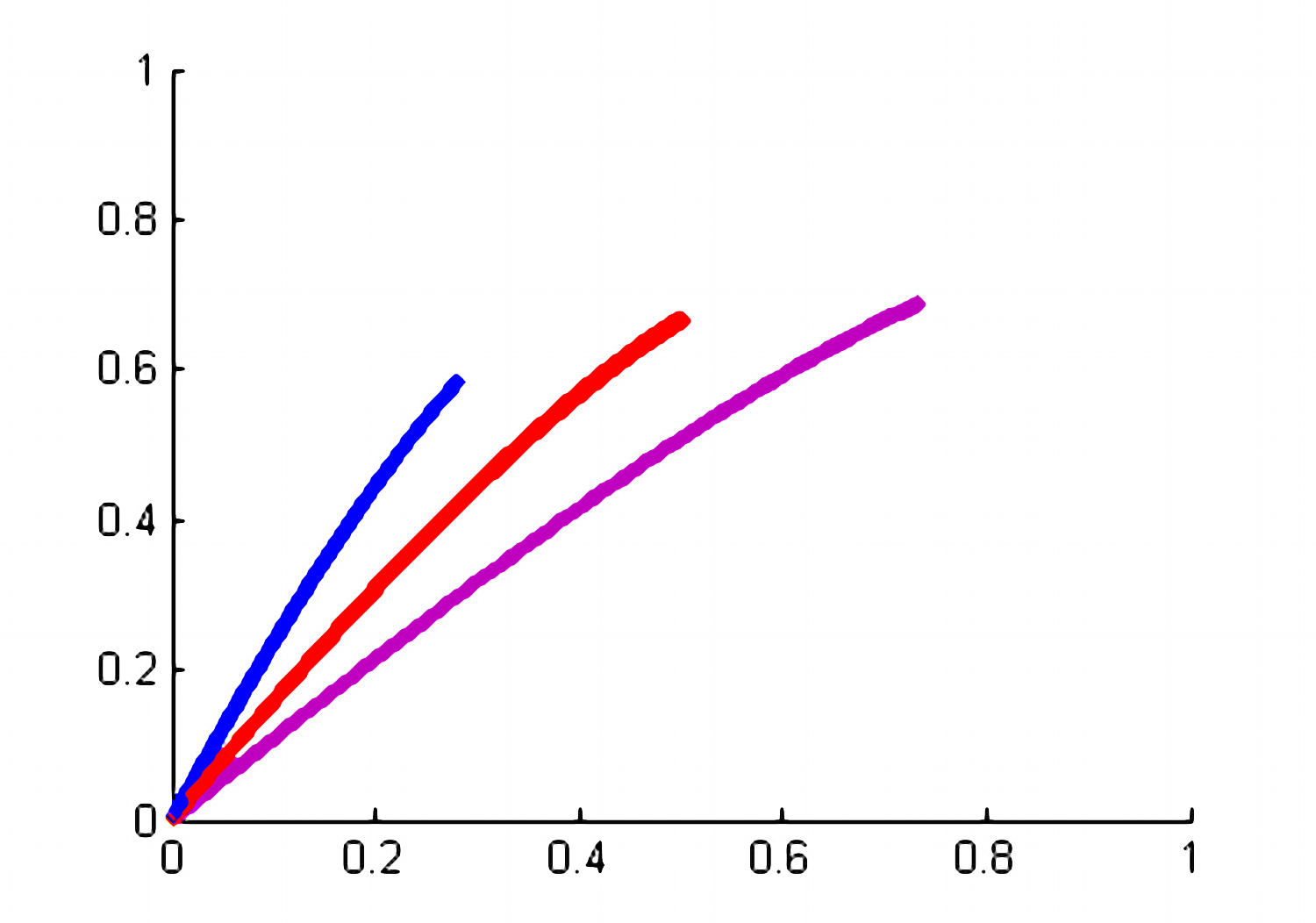}
     \caption{Three individual attractors ${\cal A}_w$ for $\Delta t = 0.005$:  left: $w=111...$; middle two (very close together): $w = 1010...$;
     right: $w=000...$.}\label{f7}
  \end{minipage}
\end{figure}

\begin{figure}[h]
\centering
  \begin{minipage}[t]{0.4\linewidth}
  \includegraphics[width=2in]{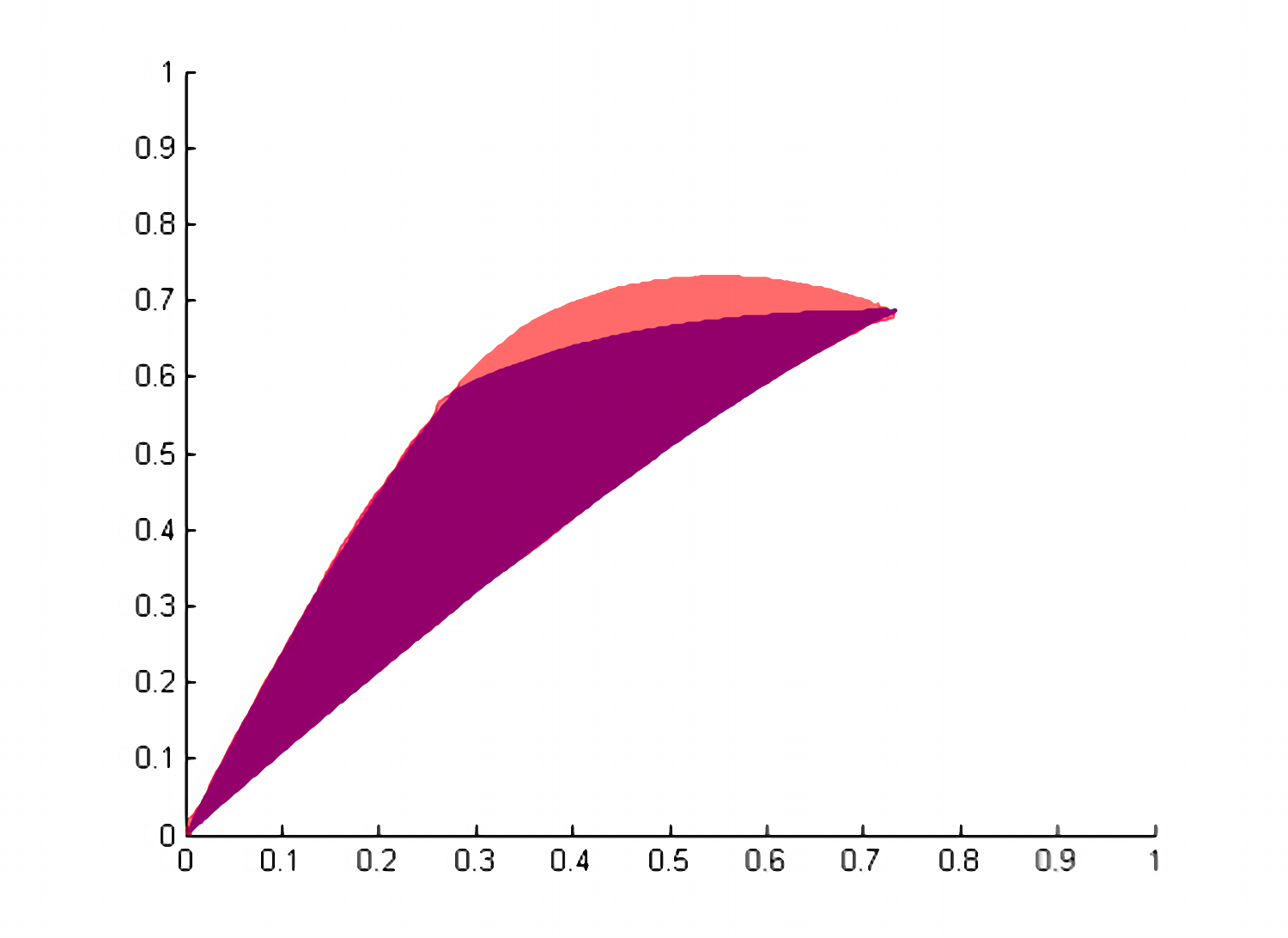}
     \caption{``Convex combination" superimposed over $K$; $\Delta t = 0.005$.}\label{f8}
     \end{minipage}%
  \hspace{0.5in}%
  \begin{minipage}[t]{0.4\linewidth}
  \includegraphics[width=2in]{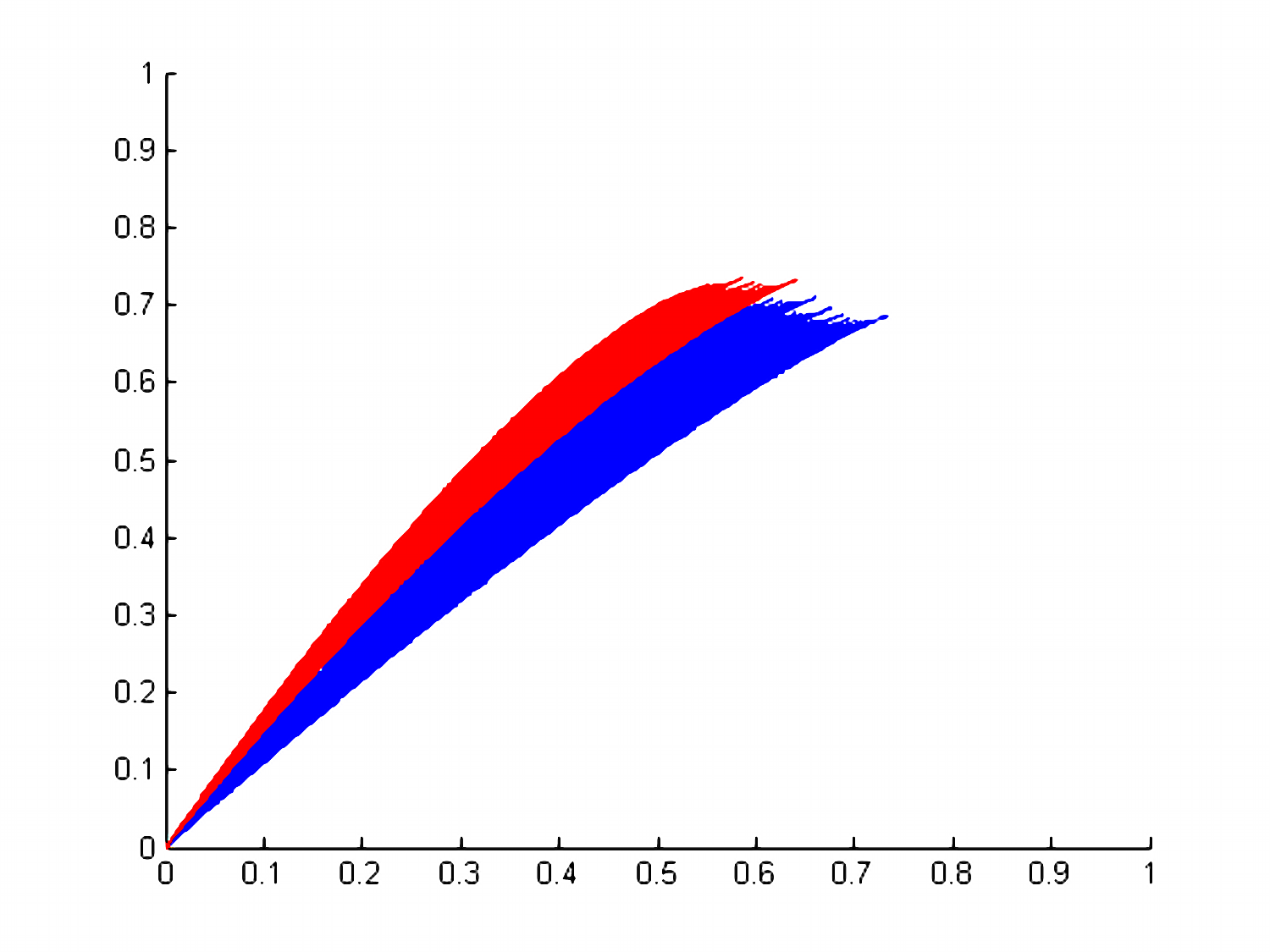}
     \caption{Golden mean, full; $\Delta t = 0.05$.}\label{f9}
     \end{minipage}
\end{figure}

We choose two sets of parameters,
$\hbox{pset}_0 = \{a = 4, \,b = 6,\,r = 1,\,m = 2\}$ and $\hbox{pset}_1 = \{a = 2, \,b = 10,\,r = 3,\,m = 2\}$, and denote the corresponding time step maps
by $S_0$ and $S_1$. These sets of parameters are
not related to any real-life situation but rather chosen to better visualize  the attractors.
The fixed point $P_2$ for $\hbox{pset}_0$ is
$(11/15,\,11/16)$ and for $\hbox{pset}_1$ it is
$(7/25,\,7/12)$.
Figures \ref{f2} through \ref{f9} show the results of numerical computation. The results depend on the
size of the time step $\Delta t$.
On figures \ref{f6} and \ref{f7}, the left line (the heteroclinic trajectory) is
the (global compact) attractor for the discrete system
$(X, S_1)$, and the right line is the attractor
of $(X, S_0)$. The two lines between them form the
individual attractor ${\cal A}_w$ corresponding to
the periodic string $w = 1010...$ (on figure \ref{f7} the two line are very close).  For our example of dynamics with choice, $\Sigma$ is the space of one-sided infinite strings of symbols $0$ and $1$. According to Theorem \ref{two}, the global compact attractor for $({\frak X}, \Sigma)$ has one slice, i.e., ${\frak M} = K\times \Sigma$. The set $K$ for $\Delta t = .05$ and for
$\Delta t = .005$ are depicted on figures \ref{f6}
and \ref{f7}, respectively.
We have also looked at the dynamical systems
corresponding to convex combinations
of the parameter sets $\hbox{pset}_0$ and
$\hbox{pset}_1$ and plotted their global attractors.
The result is different from $K$, see figure
\ref{f6} where the ``convex combination" is
superimposed onto the set $K$.

\begin{figure}[h]
\centering
\includegraphics[width=2in]{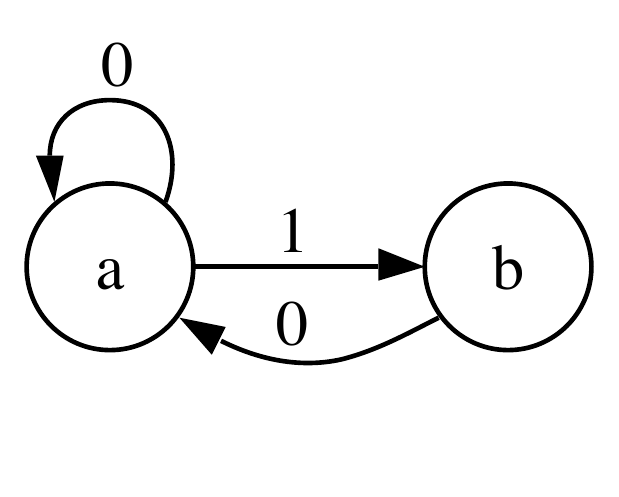}
\caption{The golden mean shift.}
\label{f10}
\end{figure}

\begin{figure}[h]
\centering
  \begin{minipage}[t]{0.4\linewidth}
     \includegraphics[width=2in]{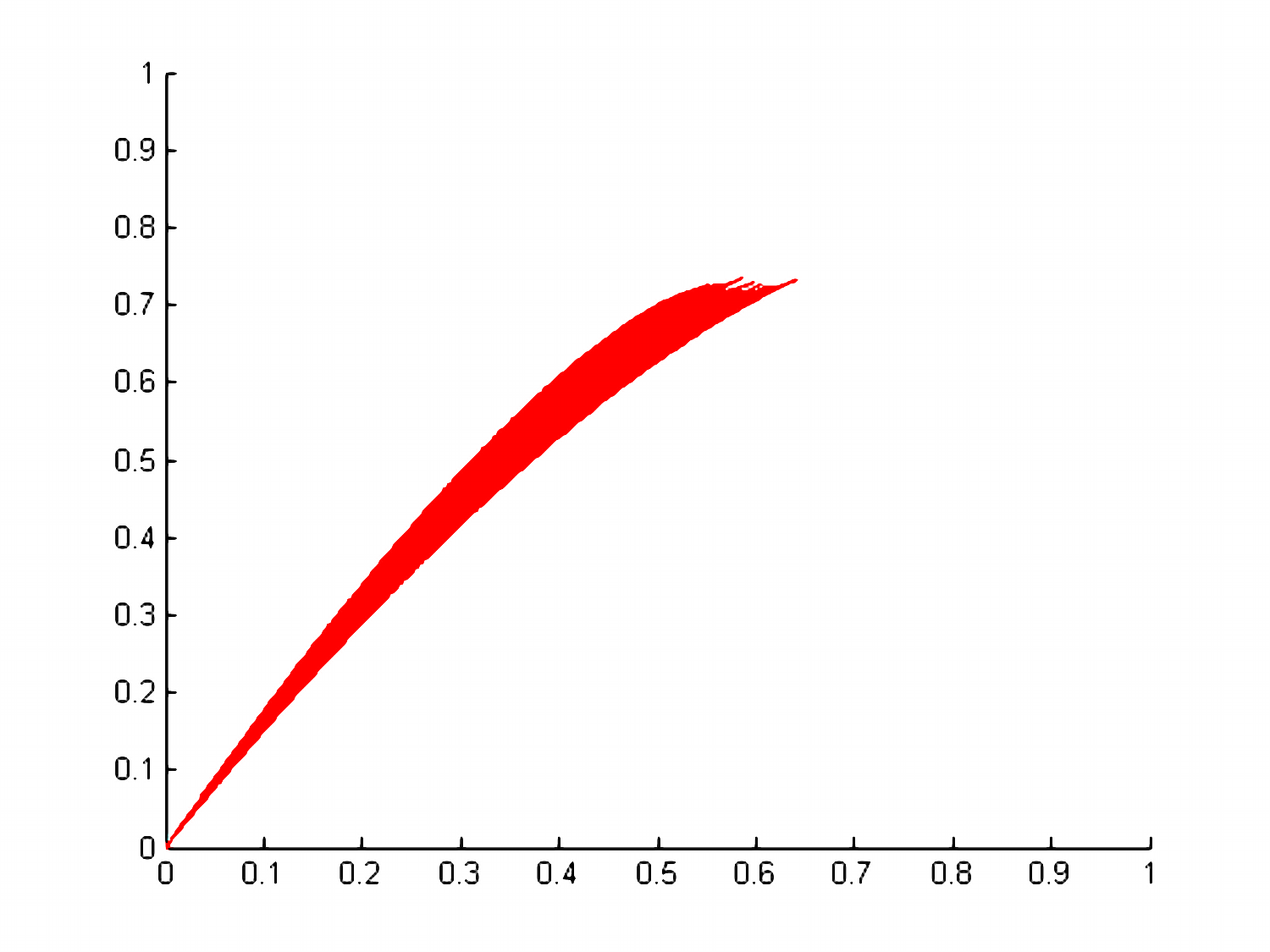}
     \caption{The red slice; $\Delta t = 0.05$.}\label{f11}
  \end{minipage}%
  \hspace{0.5in}%
  \begin{minipage}[t]{0.4\linewidth}
     \includegraphics[width=2in]{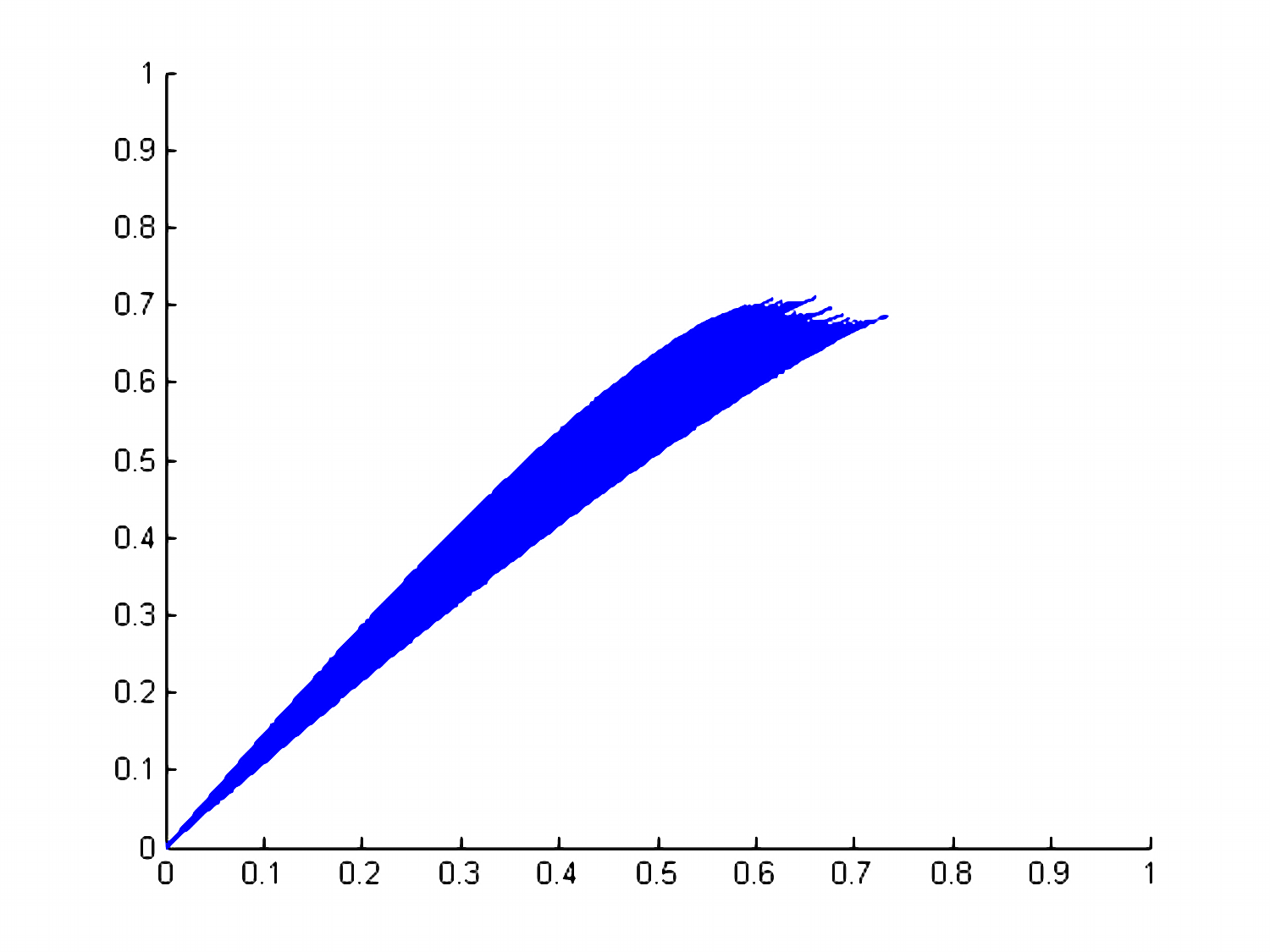}
     \caption{The blue slice; $\Delta t = 0.05$.}\label{f12}
  \end{minipage}
\end{figure}
When $\Delta t \to 0$, the upper part of the boundary of $K$
becomes smooth. Note that the limit set is not an
attractor of any system (\ref{dode}) with a fixed,
averaged set of parameters $a, b, r$ and $m$.
It would be interesting to understand whether
the limit set can be obtained as a union of
the attractors of the systems $(X, S_t)$, where
the operator $S_t$ corresponds to a certain parameter set
$\hbox{pset}_t$ for some curve connecting
$\hbox{pset}_0$ with $\hbox{pset}_1$ in the space
of parameters.

Next, we consider restricted dynamics
associated with the golden mean subshift $\Lambda$
(made of one-sided strings of $0$s and $1$s
such that each $1$ is necessarily followed by $0$).
The graph representing the golden mean shift is
shown on figure \ref{f10}.

Our analysis in section \ref{rd} shows that the global attractor of the
restricted dynamics, $({\frak X}, \Lambda)$
may have at most two different slices: one
corresponding to sequences of words ending in
$1$ (the red slice), and the other one
corresponding to sequences of words ending in
$0$ (the blue slice).
Our computation shows that the attractor
of the restricted dynamics $({\frak X}, \Lambda)$
indeed has
two slices. The slices
are shown on figures
\ref{f11} and \ref{f12}.

As point sets on the plane, the slices overlap.
Their union is plotted on figure \ref{f9}.

\end{document}